%% file: article02.tex
\input Macros.tex
\def\Vol{{\rm Vol}}
\def\vol{{\rm vol}}


\arXiv 

\anglais
\tenpoint

\TITR THE MACROSCOPICAL SOUND OF TORI|


\auteur Constantin Vernicos|


\centerline{\bf Abstract}
\medskip
\centerline{\vbox{\hsize 10cm \baselineskip=2.5ex \hautdepage Take a torus with a Riemannian metric. 
Lift the metric on its universal cover.
You get a distance which in turn yields balls. 
On these balls you can look at the Laplacian. Focus
on the spectrum for the Dirichlet or Neumann problem. 
We describe the asymptotic behaviour of the eigenvalues
as the radius of the balls goes to infinity, and characterise
the flat tori using the tools of homogenisation our conclusion being
that "Macroscopically, one can hear the shape of a flat torus". We also
show how in the two dimensional case we can recover earlier results
by D.~Burago and S.~Ivanov on the asymptotic volume.}}

\vfootnote{}{{\sevensl Key-words : }{\hautdepage Spectrum of the laplacian, Tori, Homogenisation,
Stable norm, Asymptotic volume.}}

\vfootnote{}{{\sevensl Classification : }{\hautdepage 53C24, 58C40, 74Q99}}




\partie{Introduction and claims}\label{part1}

Let $(\T^n,g)$ be a riemannian Torus, lift its metric on its universal cover and use
it to define first a distance, then the metric's balls. The first thing one can observe
is the volume of these balls as a function of their radius, 
indeed as the distance obtained arises from
a compact quotient it is equivalent to an euclidean distance hence the volume of these
balls is equivalent to the euclidean volume of an euclidean ball i.e. proportional to the radius of
the ball to the power of $n$ (the dimension of our Torus).

We are thus naturally led to wonder what happens if one
looks at the following riemanniann function on the balls ($B_g(\rho)$ is the ball of radius $\rho$)
$$
{{\rm Vol}_g\bigl(B_g(\rho)\bigr) \over \rho^n}\quad {\rm as} \quad \rho \to +\infty\ {\rm .}
$$
If it is not very
surprising that it converges to some constant for this limit can be seen
as a mean value due to the periodicity of the metric 
(see for example Pansu \cite{pansu} and a slightly different and more analytical proof in this paper section \ref{newvolas}),
it is quite remarkable that this constant, called asymptotic volume, is bounded from below by
the constant arising from the flat tori and furthermore that the case of equality caracterises
the flat tori as D.~Burago and S.~Ivanov showed in \cite{bi2}.

The study of the balls of large radii on the universal cover of tori (and more generally
of a nilmanifold) is what we call here the macroscopical geometry. Indeed in our case
the universal cover is a real vector space, where some lattices
acts by translation (in the more general case of nilmanifolds one should consider a left action). 
Should one focus on the point of this lattice endowed with the distance
arising from the torus, one gets an invariant metric on the lattice. Now if
one looks at this lattice from a galaxy far, far away, one won't be able to distinguish the lattice
from the whole universal cover. Thus it is understandable that for this observer the
distance obtained on the universal cover seems invariant by all translations
(for general nilmanifolds one gets a left invariant distance).

In the case of tori this "seen from a far away galaxy" 
distance is a norm, called the stable norm and was first
defined by Federer in Homology. It is some kind of mean value of the metric. 
This asymptotic behaviour was generalized and prooved by P.~Pansu for
all nilmanifold \cite{pansu} and precised by D.~Burago \cite{bu1} for Tori. Since then
the stable normed appeared in many other works : for surfaces and
the links with Aubry-Mather theory in D.~Massart's works, one can also find it in
the weak KAM theory of A.~Fathi. It is also worth mentionning the crucial role
it plays in the proof by D.~Burago and S.~Ivanov \cite{bi1} of the Hopf conjecture
concerning tori without conjugate points.
 Here in part \ref{part2} we show, for the case of tori, how one recovers the stable norm
using Homogenisation's tools.

There is another interesting geometric invariant attached to the balls and
linked with the volume, the spectrum of the laplacian. Indeed if one knows
the spectrum one knows the volume thanks to Weyl's asymptotic formula. Here again
one easily sees, comparing with the euclidean case, that the eigenvalues converge
to zero with a $1/\rho^2$ speed ($\rho$ being the radius). If one can expect a convergence
when rescaled, it is quite surprising that as a limit we obtain the spectrum
of an euclidean and not a finsler metric, indeed the behaviour is 
described by the following theorem which is one of the aims of this paper :

\proclaim  Theorem \theo. Let $(\T^n,g)$ be a Riemannian Torus, $B_g(\rho)$ the induced metric ball on
its universal cover and $\lambda_i\bigl(B_g(\rho)\bigr)$ the $i^{\rm th}$ eigenvalue of the Laplacian 
for the Dirichlet (resp. Neumann) problem.
\endgraf
There exists an elliptic operator $\Delta_\infty$, which is the Laplacian of some Euclidean metric
on $\R^n$, such that if $\lambda_i^\infty$ is its $i^{\rm th}$ eigenvalue for the Dirichlet (resp. Neumann) problem
on the stable's norm unit ball then \label{theo1}
$$
\lim_{\rho\to +\infty} \rho^2\lambda_i\bigl(B_g(\rho)\bigr) = \lambda_i^\infty
$$

The section \ref{proof1} is devoted to the proof of this theorem and section \ref{part3}
introduces the analytical background : Homogenisation and various convergence fairly known by
the specialist of homogenisation but adapted here to our purpose, as it becomes
after some re-spelling of the problem made in section \ref{respell}.
From this theorem
we can deduce the following one which is some kind of ``{\`a} la'' Burago-Ivanov macroscopical
rigidity and which inspired the abstract :

\proclaim Theorem \theo. Let $(\T^n,g)$ be a Riemannian Torus, $B_g(\rho)$ the induced metric ball on
its universal cover and $\lambda_1\bigl(B_g(\rho)\bigr)$ the first eigenvalue of the Laplacian 
for the Dirichlet problem. Then
\item{(a)} $\displaystyle{\lim_{\rho\to +\infty} \rho^2\lambda_1\bigl(B_g(\rho)\bigr) = \lambda_1^\infty \leq \lambda_{e,n}}$
\item{(b)} equality holds if, and only if, the torus is flat
\endgraf
\noindent where $\lambda_{e,n}$ is the first eigenvalue of the Euclidean Laplacian on the Euclidean
 unit ball.\label{theo2}

The proof which is done in section \ref{proof2}, involves
some kind of transplantation for the inequality mixed with $\Gamma$-convergence
for the equality. 
For a better understanding of
what happens we briefly give some informations related to the $\Gamma$-convergence and 
adapt it to our purpose in section \ref{topo2}, following the general ideas
of K.~Kuwae and T.~Shioya in \cite{kushi} (who in turn generalized U.~Mosco's paper \cite{mosco}),
this section being completed by the proof of section \ref{lastbutnotleast}. 

As the macroscopical spectrum involved rises from an euclidean metric, we can use
 the Faber-Krahn inequality to obtain a new inequality
regarding the asymptotic volume, this is done in  section \ref{monvolas} :

\proclaim Proposition \theo.
Let $(\T^n,g)$ be a Riemannian Torus, $B_g(\rho)$ the geodesic balls of radius $\rho$
centred on a fixed point and ${\rm Vol}_g\bigl(B_g(\rho)\bigr)$ their Riemannian volume 
induced on the universal cover, writing \label{burivame}
$$
{\rm Asvol}(g) =\lim_{\rho\to \infty} \frac{{\rm Vol}_g\bigl(B_g(\rho)\bigr)}{\rho^n}
$$
then
\item{1.} $ {\rm Asvol}(g) \geq \frac{\displaystyle{{\rm Vol}_g(\T^n)}}{\displaystyle{{\rm Vol}_{{\rm Al}}(\T^n)}}\ \omega_n $
\item{2.} In case of equality, the torus is flat.
\endgraf
\noindent Where $\omega_n$ is the unit Euclidean ball's Euclidean volume,
and ${\rm Vol}_{{\rm Al}}(\T^n)$ is the volume of the albanese Torus.

 A surprising facts arises
because this new inequality involves a constant which happened to be at the
heart of the isosystolic inequality of two dimensional tori (see J.~Lafontaine \ref{jaklaf}),
hence we obtain an alternate proof of the asymptotic volume's lower 
boundedness in dimension two,

\proclaim Corollary.
Let $(\T^2,g)$ be a $2$-dimensional torus then :
\item{1.} $ {\rm Asvol}(g) \geq \pi $
\item{2.} In case of equality, the torus is flat.

It is worth mentionning that the case of equality in the previous
claims relies on theorem \ref{sneuc}, which states that the stable norm
coincides with the Albanese metric if and only if the torus is flat, and
whose proof do not depend on the work of D.~Burago and 
S.~Ivanov \cite{bi2}, thus we actually get an alternate proof of their theorem in
the $2$-dimensional case.

We also give a kind
of generalised Faber-Krahn inequality for normed finite dimensional vector spaces, 
which implies that we cannot distinguish  the euclidean's ones among them using
the first generalised eigenvalue of the Dirichlet laplacian (see lemma \ref{lemfk} and its corollary):

\proclaim Theorem {\theo.} (Faber-Krahn inequality for norms). 
Let $D$ be a domain of $\R^n$, with the norm $\n{\cdot}$ and a measure $\mu$ invariant
by translation. Let $D^*$ be the norm's ball with same measure as $D$, then
$$
\lambda_1\bigl(D^*, \n{\cdot}\bigr) \leq %
\lambda_1\bigl(D, \n{\cdot}\bigr)
$$

We finally explain in section \ref{lgtime} how is our work related to works focused on the long time
asymptotics of the heat kernel (see \cite{kosu}, \cite{duzu}, \cite{gconv})
and finally in section \ref{grnilm} we state how theorem \ref{theo1} transposes to all 
graded nilmanifolds (subject which should be widely extended in a forthcoming article).

I would like to thanks my advisor G{\'e}rard Besson for his good advice and
Y. Colin~de~Verdi{\`e}re who lead me to the $\Gamma$-convergence which happened to be the
adapted tool, as I hope to convince the reader, for the macroscopical geometry.
\partie{Stable Norm and Homogenisation}\label{part2}
In this section we show how the stable norm, the Gromov-Haussdorff convergence
and the $\Gamma$-convergence of the homogenisation theory are linked and finish
by re-spelling our goal. In what follows, $B_g(\rho)$ will be the metric ball
of radius $\rho$ on the universal cover of a torus with the lifted metric. We
first begin by two definitions.

\subpartie{Convergences}
We recall the definition of $\Gamma$-convergence in a metric space :

\proclaim Definition \theo. Let $(X,d)$ be a metrics space. We say that a sequence of 
function $(F_j)$ from $X$ to $\overline\R$, $\Gamma$-converges to a function $F:X\to\overline\R$ 
if and only if for all $x\in X$  we have \label{gamacv1}
\item{1.} for all converging sequence $(x_j)$ to $x$
$$
F(x)\leq \liminf_{j\to \infty} F_j(x_j);
$$
\item{2.} there exists a sequence $(x_j)$ converging to $x$ such that
$$
F(x)= \lim_{j\to \infty} F_j(x_j).
$$

We now introduce the Gromov-Hausdorff measured convergence in the space ${\cal M}$ of compact metric 
and measured spaces $(X,d,m)$ modulo isometries. 
First if $X$ and $Y$ are in ${\cal M}$ then an application
$\phi:X\to Y$ is called an {\sl $\epsilon$-Hausdorff approximation} if and only if we have
\item{(1.)} the $\epsilon$-neighbourhood of $\phi(X)$ in $Y$ is $Y$;
\item{(2.)} for all $x,y\in X$ we have
$$
\Bigl| d(x,y)-d\bigl(\phi(x),\phi(y)\bigr) \Bigr| \leq \epsilon.
$$

We write $C^0(X)$ for the space of continues functions from $X$ to $\R$ and ${\cal A}$ will
be a partially directed space.

\proclaim Definition \theo. We say that a net $(X_\alpha,d_\alpha,m_\alpha)_{\alpha\in {\cal A}}$ of spaces in ${\cal M}$
converges to $(X,d,m)$ for the Gromov-Hausdorff measured topology if, and only if there exists a net
of positive real numbers $(\varepsilon_\alpha)_{\alpha\in {\cal A}}$ decreasing to $0$ and $m_\alpha$ measurable $\varepsilon_\alpha$-Hausdorff 
approximations $f_\alpha:X_\alpha\to X$ such that $(f_\alpha)_*(m_\alpha)$ converges vaguely to $m$ i.e. 
$$
\int_{X_\alpha} u\circ f_\alpha~dm_\alpha \to \int_X u~dm \quad \forall u \in C^0(X).
$$

\subpartie{The stable norm}\label{sectionnorm}
Let $(\T^n,g)$ a Riemannian torus. We will call {\sl rescaled metrics}
the metrics $g_\rho=(1/\rho^2)(\delta_\rho)^*g$ and their lifts on the universal cover. We will
also write $\delta_\rho$ for the homothetie of scale $\rho$.

In the $80$'s P.~Pansu showed that the distance induced on $\R^n$ as a universal
cover of a torus, behaved asymptotically like the distance induced by a norm. In the
$90$'s D.~Burago showed a similar result for periodic metrics on $\R^n$. It is that norm
which is called the {\sl stable norm}. To be more precise let us write $f_1(x)=d_g(0,x)$ the
distance from the origin to $x$ and $f_\rho(x)=d_g(0,\delta_\rho(x))/\rho$, then P.~Pansu's result says
that there exist a norm $||.||_\infty$ such that for all $x\in \R^n$
$$
\lim_{\rho\to \infty} f_\rho(x) = ||x||_\infty
$$ 
and Burago's says that there exists a constant $C$ such that for all $x\in \R^n$
$$
\bigl|f_\rho(x)-||x||_\infty\bigr| \leq {C\over \rho}
$$
in other words, Pansu's results is a simple convergence and Burago's is a uniform convergence result.

There is another proof of the simple convergence of the sequence $(f_\rho)$ as $\rho$ goes to infinity,
using homogenisation's tools.

\proclaim Theorem \theo. Let $\tilde g$ the induced metric on $\R^n$ as a universal 
cover of a Riemannian Torus $(\T^n,g)$,
then there exists a norm $||.||_\infty$ such that
\item{1.} for every bounded open $I\subset\R$ the sequence of functional on $W^{1,2}(I;\R^n)$
$$
E_{\rho}(u)=\int_I \tilde g_{(\delta_\rho u(t))}\bigl( u'(t),u'(t)\bigr)dt
$$
$\Gamma$-converge for the $L^2$ norm toward the functional
$$
E_\infty(u)=\int_I \bign{u'(t)}_\infty^2 dt;
$$
\item{2.} the norm satisfies
$$
\n{\xi}_\infty= \lim_{t\to +\infty } \inf \biggl\{ {1 \over t} \int_0^t \tilde{g}_{(u+\xi \tau)} (u'+\xi,u'+\xi)~d\tau :%
u \in W^{1,2}_0\bigl(\mathopen]0,t\mathclose[;\R^n\bigr) \biggr\}. \numeq
$$
\endgraf
Furthermore if $f_\rho(x)=d_g(0,\rho x)/\rho$ then for all $x \in \R^n$
$$
\lim_{\rho\to +\infty}f_\rho(x)=\n{x}_\infty
$$

\proof  We use the Proposition 16.1 page 142 of A.~Braides and A.~Desfranceschi \cite{brdef}.
It gives us the $\Gamma$-convergence of the sequence of functional $(E_\rho)$ toward a functional $E_\infty$
such that
$$
E_\infty(u)=\int_I \varphi\bigl(u'(t)\bigr) dt
$$
with $\varphi$ convex and satisfying the asymptotic formula (1). It remains to show
that $\varphi$ is the square of a norm.

{\bf Homogeneity :} Using the asymptotic formula (1) we easily get $\varphi(0)=0$ and by a change
of variables $\varphi(\lambda x)=\lambda^2\varphi(x)$.

{\bf Separation :} Let us point that
\item{1.} The minimum of the energy of a path between $0$ and $t\xi$ in an Euclidean space
is attained for the straight line. Thus if we put into the asymptotic formula (1) an Euclidean metric,
we get the same metric.
\item{2.} Let $g$ and $h$ be two metrics such that for all $s$ and $\xi$
$$
g_s(\xi,\xi)\leq h_s(\xi,\xi) 
$$
then for all $u\in  W^{1,2}_0(\mathopen]0,t\mathclose[;\R^n)$ we get
$$
 \frac{1}{t} \int_0^t g_{(u+\xi \tau)} (u'+\xi,u'+\xi)~d\tau \leq \frac{1}{t} \int_0^t h_{(u+\xi \tau)} (u'+\xi,u'+\xi)~d\tau
$$
thus taking the infimum for $u$ and taking the limit as $t$ goes to infinity we get
$$\displaylines{%
\qquad \lim_{t\to +\infty } \inf \biggl\{ \frac{1}{t} \int_0^t g_{(u+\xi \tau)} (u'+\xi,u'+\xi)~d\tau :%
u \in W^{1,2}_0\bigl(\mathopen]0,t\mathclose[;\R^n\bigr) \biggr\} \leq \hfill\cr
\hfill \lim_{t\to +\infty } \inf \biggl\{ \frac{1}{t} \int_0^t h_{(u+\xi \tau)} (u'+\xi,u'+\xi)~d\tau :%
u \in W^{1,2}_0\bigl(\mathopen]0,t\mathclose[;\R^n\bigr) \biggr\}. \qquad \nume}
$$

Now let us also remark that $g$ being periodic, there exists two strictly positive 
constants $\alpha$ and $\beta$ such that
$$
\alpha|\xi|^2\leq g_s(\xi,\xi)\leq\beta|\xi|^2
$$
now applying the three remarks we get
$$
\alpha|\xi|^2\leq \varphi(\xi)\leq\beta|\xi|^2
$$
thus $\varphi(\xi)=0$ if and only if $\xi=0$.

{\bf Triangle inequality : } First note that
$$\{\xi \in \R^n \mid \varphi(\xi)\leq1\}=\{\xi \in \R^n \mid \sqrt{\varphi(\xi)}\leq1\}=S_n$$ 
It follows that if $\sqrt{\varphi}(\xi)=1=\varphi(\xi)$
and $\sqrt{\varphi}(\eta)=1=\varphi(\eta)$ then for all $0\leq\lambda\leq1$ by the convexity of $\varphi$ 
$$
\varphi(\lambda\xi + (1-\lambda) \eta) \leq \lambda\varphi(\xi)+(1-\lambda)\varphi(\eta) =1
$$
so
$$
\sqrt{\varphi}(\lambda\xi + (1-\lambda) \eta) \leq1.
$$
Thus for all non null $x$, $y$ 
$$
\sqrt{\varphi}\Bigl(\lambda \frac{x}{\sqrt{\varphi}(x)} + (1-\lambda)\frac{y}{\sqrt{\varphi}(y)}  \Bigr) \leq 1
$$
now taking $\lambda=\sqrt{\varphi}(x) / \bigl(\sqrt{\varphi}(x)+\sqrt{\varphi}(y)\bigr)$ and using
$\sqrt{\varphi}$ homogeneity we finally get the triangle inequality and we are able
to conclude that $\n{\cdot}_\infty=\sqrt{\varphi}(\cdot)$ is a norm.

The final assertion comes from the fact that $\n{\xi}_\infty^2$ is the limit of the energies' infimum
along the paths between $0$ and $\xi$ for the rescaled metrics $(1/t^2)(\delta_t^*)g$, which are attained
along the geodesics. \qed

This theorem easily induces the following assertion.
\proclaim Corollary \theo. For all $x$ and $y\in \R^n$ we have \label{cor:normes}
$$
\lim_{\rho\to +\infty}{d_g(\rho x,\rho y) \over \rho} =\n{x-y}_\infty
$$

From now on we will write $d_\rho(x,y)=d_g(\rho x,\rho y)/\rho$, and we are now going to see
what can be deduced for the balls $B_g(\rho)$ in terms of Gromov-Haussdorff convergence.

\subpartie{Gromov-Hausdorff convergence of metric balls}\label{newvolas}
We will write $\mu_g$ (resp. $\mu_\rho$) the measure induced by $\tilde g$ (resp. $g_\rho$). 
$\mu_\infty$ will be the measure of Lebesgue such that for a fundamental domain $D_f$
we have $\mu_\infty(D_f)=\mu_g(D_f)$. Finally let 
$$B_\rho(R)=\bigl\{ x\in\R^n \bigm| d_\rho(0,x)\leq R  \bigr\}={1\over \rho }B_g(R\cdot\rho),$$
and 
$$B_\infty(R)= \bigl\{ x\in\R^n \bigm| \n{x}_\infty\leq R  \bigr\}.$$

\proclaim Theorem \theo. The net of measured metric spaces $\bigl(B_\rho(1),d_\rho,\mu_\rho\bigr)$ converges
in the Gromov-Hausdorff measured topology to $\bigl(B_\infty(1),\n{\cdot}_\infty,\mu_\infty\bigr)$  
as $\rho$ goes to infinity.\label{lesfonctionsalg}

\proof Let us choose an $\epsilon>0$. We first show that the identity is an $\epsilon$-approximation
if $\rho$ is large enough. It suffice to show that there is a finite family of points $(x_1,\ldots,x_N)$
such that its $\epsilon$-neighbourhood in $\bigl(B_\infty(1),d_\infty\bigr)$ and, for $\rho$ large enough, 
in $\bigl(B_\rho(1),d_\rho\bigr)$ is respectively $B_\infty(1)$ and $B_\rho(1)$ and such that for all $i,j =1,\dots,N$
we have
$$
\bigl| \n{x_i-x_j}_\infty -d_\rho(x_i,x_j) \bigr| \leq \epsilon
$$

Let $r>0$ and let $(\gamma_1,\ldots,\gamma_N)$ be all the images of $0$ by the action of $\Z^n$, such that
for $i=1,\ldots,N$, $\gamma_i \in B_\infty(r)$. Then we take for $i=1,\ldots,N$, $x_i=\gamma_i/r$. Let us remark that
for $\rho$ large enough these points will all be in $B_\rho(1)$.

Now let us point out that, because of the invariance by the $\Z^n$ action, there are two constants
$\alpha$ and $\beta$ such that for all $x$ and $y\in \R^n$ we have
$$
\alpha\n{x-y}_\infty \leq d_g(x,y) \leq \beta \n{x-y}_\infty ;
$$
thus for every $x\in B_\infty(1)$ take the closest point $x_i$ 
(thus $\gamma_i$ is the closest point of $\Z^n\cdot0$ from $rx$)  then there is a constant $C$ (the diameter
of the fundamental domain) such that 
$$
\n{x-x_i}_\infty\leq {1 \over \alpha r}d_g(rx,\gamma_i)\leq {1\over \alpha r}C
$$
we also get
$$
d_\rho(x,x_i) \leq {\beta\over \alpha r}C
$$
thus, for $r$ large enough $(x_1,\ldots,x_N)$ is an $\epsilon$-neighbourhood of $\bigl(B_\infty(1),\n{\cdot}_\infty\bigr)$. Furthermore
if $\rho$ is large enough it is also an $\epsilon$-neighbourhood of $\bigl(B_\rho(1),d_\rho\bigr)$ and by corollary 6
$$
\bigl| \n{x_i-x_j}_\infty -d_\rho(x_i,x_j) \bigr| \leq \epsilon.
$$

Now let us take a continuous function from $B_\infty(1)$ to $\R$. Let $z_1,\ldots,z_k$ and $\zeta_1,\ldots,\zeta_l$
in the orbit of $0$ by the $\Z^n$ action such that $\zeta_j+D\cap B_\infty(\rho)\not= \emptyset$ for $j=1,\ldots l$ and
$$
\bigcup_i z_i+D_f \subset B_\infty(\rho) \subset \bigcup_k \zeta_k+D_f
$$
(where we took all $z_i$ such that $z_i+D_f \subset B_\infty(\rho)$) then we get
$$\displaylines{%
\qquad \sum_i\inf_{\rho x \in z_i+D_f} f(x)~\mu_\infty(D_f) \leq \hfill \cr
\hfill \int_{B_\infty(\rho)} f(x/\rho)~d\mu_g(x) \leq \sum_j\sup_{\rho x \in (\zeta_j+D_f)\cap B_\infty(\rho)} f(x)~\mu_\infty(D_f) \qquad }
$$
now dividing by $\rho^n$ we find
$$\displaylines{%
\qquad \sum_i\inf_{x \in {1\over \rho}(z_i+D_f)} f(x)~\mu_\infty\bigl((1/\rho) D_f\bigr)\leq \hfill \cr
\hfill \int_{B_\infty(1)} f ~d\mu_\rho(x) \leq \sum_j\sup_{x \in {1\over \rho}(\zeta_j+D_f)\cap B_\infty(1)} f(x)~\mu_\infty\bigl((1/\rho) D_f\bigr) \qquad
}$$
The middle term is surrounded by two sums of Riemann, which converges to $\int_{B_\infty(1)}f~d\mu_\infty$,
thus it also converges. To conclude, notice that the net of characteristic
function $\chi_{B_\rho(1)}$ converges simply to $\chi_{B_\infty(1)}$ inside of $B_\infty(1)$. \qed

\subpartie{What shall we finally study ?} \label{respell}

As we said we are now going to focus on the spectrum of the balls $B_g(\rho)$. As we allready
mentioned we know that the eigenvalues are converging to zero with a $1/\rho^2$ speed. Hence we
want to find a precise equivalent.

For this let introduce $\Delta_\rho$ the laplacian associated to the rescaled metrics $g_\rho=1/\rho^2(\delta_\rho)^*g$,
and for any function $f$ from $B_g(\rho)$ to $\R$ lets associate a function $f_\rho$ on
$B_\rho(1)$ by $f_\rho(x)=f(\rho\cdot x)$. Then it is an easy calculation to see that for any $x\in B_\rho(1)$ :
$$
\rho^2\bigl(\Delta f\bigr)(\rho\cdot x)= \bigl(\Delta_\rho f_\rho\bigr)(x)
$$
hence the eigenvalues of $\Delta_\rho$ on $B_\rho(1)$ are exactly the eigenvalues of $\Delta$ on $B_g(\rho)$
multiplied by $\rho^2$ and
our problems becomes the study of the spectrum of the laplacian
$\Delta_\rho$ on $B_\rho(1)$. In the light of what precedes we would like to show that
there is some operator $\Delta_\infty$ acting on $B_\infty(1)$ such that, in some sense, the net
of laplacian $(\Delta_\rho)$ converges towards $\Delta_\infty$ such that the spectra also
converge to the spectrum of $\Delta_\infty$. The following section aims at giving a
precise meaning to this.

\partie{Convergence of spectral nets}\label{part3}

This section adapts to our purpose some notion of convergences
well known for a fixed Hilbert space.

\subpartie{Convergence on a net of Hilbert spaces}
Let $(X_\alpha,d_\alpha,m_\alpha)_{\alpha\in {\cal A}}$, where ${\cal A}$ is a partially ordered set, 
be a net of compact measured metric spaces converging in
the Gromov-Hausdorff measured topology to $(X_\infty,d_\infty,m_\infty)$. We will write $L_\alpha^2=L^2(X_\alpha,m_\alpha)$
(resp. $L^2_\infty(X_\infty,m_\infty)$) for the square integrable function spaces. Their respective scalar product
will be $\langle\cdot,\cdot\rangle_\alpha$ (resp.  $\langle\cdot,\cdot\rangle_\infty$) and $\n{\cdot}_\alpha$ (resp. $\n{\cdot}_\infty$). 

Furthermore we suppose that in every  $L_\alpha^2$ the continuous functions form a dense subset.

\proclaim Definition \theo. We say that a net $(u_\alpha)_{\alpha \in {\cal A}}$ of functions $u_\alpha\in L^2_\alpha$ strongly
converges to $u\in L^2_\infty$ if there exists a net $(v_\beta)_{\beta\in {\cal B}} \subset C^0(X_\infty)$ converging to $u$ in $L^2_\infty$
such that
$$
\lim_\beta \limsup_\alpha \n{f^*_\alpha v_\beta -u_\alpha}_\alpha =0;
$$
where $(f_\alpha)$ is the net of Hausdorff approximations. We will also talk of strong convergence in 
${\cal L}^2$.

\proclaim Definition \theo. We say that a net $(u_\alpha)_{\alpha \in {\cal A}}$ of functions $u_\alpha\in L^2_\alpha$ weakly
converges to $u\in L^2_\infty$ if and only if for every net $(v_\alpha)_{\alpha\in {\cal A}}$ strongly converging to $v\in L^2_\infty$
we have 
$$
\lim_\alpha \langle u_\alpha,v_\alpha\rangle_\alpha =\langle u,v\rangle_\infty \numeq
$$
We will also talk of weak convergence in ${\cal L}^2$.

The following lemmas justifies those two definitions

\proclaim Lemma \theo.
Let $(u_\alpha)_{\alpha \in {\cal A}}$ be a net of functions $u_\alpha\in L^2_\alpha$. If $(\n{u_\alpha}_\alpha)$ is
uniformly bounded, then there exists a weakly converging subnet.\label{subwcvnet}

\proof Let $(\phi_k)_{k\in \N}$ be a complete orthonormal basis of $L^2_\infty$. Using the density of continuous
functions in $L^2_\infty$, for each $k$ we can retrieve a net of continuous functions $(\varphi_{k,\beta})_{\beta\in {\cal B}}$
strongly converging to $\phi_k$ in $L^2_\infty$. Replacing by a subnet of ${\cal A}$ and ${\cal B}$ if
necessarily, we can assume that the following limit exists
$$
\lim_\beta \lim_\alpha \langle u_\alpha,f_\alpha^*\varphi_{1,\beta}\rangle_\alpha = a_1 \in \overline{\R}
$$
and from the uniform bound hypothesis it follows that $a_1 \in \R$. Repeating this procedure
we can assume that for every $k\in \N$ the following limit exists
 $$
\lim_\beta \lim_\alpha \langle u_\alpha,f_\alpha^*\varphi_{k,\beta}\rangle_\alpha = a_k \in \R {\rm .}
$$
Let us fix an integer $N$. For any $\epsilon>0$ there is a $\beta_\epsilon\in{\cal B}$ such that 
$$
\bigl|\langle\varphi_{k,\beta},\varphi_{l,\beta}\rangle_\infty-\delta_{kl} \bigr| < \epsilon
$$
for any $\beta\geq\beta_\epsilon$ and $k,l=1,\ldots,N$. Moreover for any $\beta\geq\beta_\epsilon$ there is an $\alpha_{\epsilon,\beta}\in {\cal A}$ such that
$$
\bigl|\langle f_\alpha^*\varphi_{k,\beta},f_\alpha^*\varphi_{l,\beta}\rangle_\alpha - \delta_{kl}\bigr| < 2\epsilon
$$
for any $\alpha\geq\alpha_{\epsilon,\beta}$ and $k,l=1,\ldots,N$. Let $L_{\alpha,\beta}={\rm Vect}\{f_\alpha^*\varphi_{k,\beta} \mid k=1,\ldots,N \}$ and 
$P_{\alpha,\beta}:L^2_\alpha\to L_{\alpha,\beta}$ be the projection to the linear subspace $L_{\alpha,\beta}\subset L_\alpha^2$ we have
$$
\Biggl|\sum_{k=1}^N \bigl|\langle u_\alpha,f_\alpha^*\varphi_{k,\beta}\rangle_\alpha\bigr|^2 - \n{P_{\alpha,\beta}u_\alpha}_\alpha^2\Biggr| \leq \theta_N(\epsilon)
$$
for every $\alpha \geq\alpha_{\epsilon,\beta}$ and $\beta \geq  \beta_\epsilon$, where $\theta_N$ is a function depending only of $N$ such that
$\lim_{\epsilon\to0}\theta_N(\epsilon)=0$. This implies for every $N$
$$\eqalign{
 \sum_{k=1}^N|a_k|^2 &= \lim_\beta \lim_\alpha \sum_{k=1}^N \bigl|\langle u_\alpha,f_\alpha^*\varphi_{k,\beta}\rangle_\alpha\bigr|^2 = \lim_\beta \lim_\alpha \n{P_{\alpha,\beta}u_\alpha}_\alpha^2 \cr
                &\leq \limsup_\alpha \n{u_\alpha}_\alpha^2 < \infty\cr}
$$
thus
$$
u=\sum_{k=1}^N a_k\phi_k \in L^2_\infty {\rm .} 
$$
We shall prove that some subnet of $(u_\alpha)_\alpha$ weakly converges to $u$. Take any $v\in L^2_\infty$
and set $b_k=\langle v,\phi_k\rangle_\infty$. By the properties of the strong convergence it is enough to show $(3)$
for a well chosen net. Let $v_\beta^N= \sum_{k=1}^N b_k \varphi_{k,\beta}$. By construction $v_\beta^N\in {\cal C}^0$ and 
$\lim_{N\to\infty} \lim_\beta v_\beta^N = v$ strongly. We have
$$
\lim_\beta \lim_\alpha \langle u_\alpha,f_\alpha^*v_\beta^N\rangle_\alpha = \lim_\beta \lim_\alpha  \sum_{k=1}^N b_k \langle u_\alpha,f_\alpha^*\varphi_{k,\beta}\rangle_\alpha = \sum_{k=1}^N a_kb_k
$$
which tends to $\langle u,v\rangle_\infty$ as $N\to \infty$. Thus, there exists a net of integers $(N_\beta)_\beta$
tending to $+\infty$ such that  $v_\beta^{N_\beta}$ strongly converges to $v$ and
$$
\lim_\beta \lim_\alpha \langle u_\alpha,f_\alpha^*v_\beta^{N_\beta}\rangle_\alpha = \langle u,v\rangle_\infty {\rm .}
$$
\qed
\proclaim Lemma \theo.
Let $(u_\alpha)_{\alpha \in {\cal A}}$ be a weakly converging net to $u\in L^2_\infty$. Then
$$
\sup_\alpha \n{u_\alpha}_\alpha < \infty \quad {\it and }\quad \n{u}_\infty \leq \liminf_\alpha \n{u_\alpha}_\alpha.
$$
Furthermore, the net strongly converges if and only if
$$
\n{u}_\infty = \lim_\alpha \n{u_\alpha}_\alpha.\label{wcvnet}
$$

\proof Let suppose that the $(u_\alpha)$ is weakly converging and $\sup_\alpha \n{u_\alpha}_\alpha = + \infty$. We can extract
a sequence such that $\n{u_{\alpha_k}}_{\alpha_k}>k$. Setting
$$
v_k = {1\over k}{u_{\alpha_k}\over \n{u_{\alpha_k}}_{\alpha_k}}
$$
one has $\n{v_k}_{\alpha_k}=1/k\to 0$ thus $v_k$ strongly converges to $0$, which implies
$$
\langle u_{\alpha_k},v_k\rangle_{\alpha_k}\to \langle u,0\rangle_\infty = 0
$$
but we also have
$$
\langle u_{\alpha_k},v_k\rangle_{\alpha_k} = {1\over k} \n{u_{\alpha_k}}_{\alpha_k} \geq 1
$$
this is a contradiction and thus we obtain $\sup_\alpha \n{u_\alpha}_\alpha<\infty$.

Let $(w_\alpha)_\alpha$ be a strongly converging net to $u$, then
$$\eqalign{
0  &\leq \liminf_\alpha \n{u_\alpha-w_\alpha}_\alpha^2 \cr
    &= \liminf_\alpha \bigl(\n{u_\alpha}_\alpha^2+\n{w_\alpha}_\alpha^2-2\langle u_\alpha,w_\alpha\rangle_\alpha \bigr) \cr
    &= \liminf_\alpha \n{u_\alpha}_\alpha^2 - \n{u}^2_\infty \cr}
$$
The final claim comes from the properties of the strong convergence and the following equality 
$$
\n{u_\alpha-w_\alpha}_\alpha^2=\n{u_\alpha}_\alpha^2 + \n{w_\alpha}_\alpha^2 -2\langle u_\alpha,w_\alpha\rangle_\alpha {\rm .}
$$
\qed
\subpartie{Convergence of bounded operators}
Let ${\cal L}(L^2_\alpha)$ bet the set of linear bounded operators acting on $L^2_\alpha$ and
$\n{\cdot}_{{\cal L}_\alpha}$ their norm (for $\alpha\in{\cal A}\cup{\infty}$). 
Let $B_\infty\in {\cal L}(L^2_\infty)$ and $B_\alpha\in {\cal L}(L^2_\alpha)$ for every $\alpha\in {\cal A}$.
\proclaim Theorem and Definition \theo. Let $u,v \in L^2_\infty$ and $(u_\alpha)_{\alpha\in{\cal A}}$, $(v_\alpha)_{\alpha\in{\cal A}}$ two nets
such that $u_\alpha, v_\alpha\in L^2_\alpha$. We say that the net of operators $(B_\alpha)_{\alpha\in{\cal A}}$ strongly 
(resp. weakly, compactly) converges to $B$ if $B_\alpha u_\alpha \to Bu$ strongly (resp. weakly, strongly) for
every net $(u_\alpha)$ strongly (resp. weakly, weakly) converging to $u \iff$
$$
\lim_\alpha \langle B_\alpha u_\alpha,v_\alpha\rangle_\alpha = \langle Bu,v\rangle_\infty \numeq
$$
for every $(u_\alpha)$, $(v_\alpha)$, $u$ and $v$ such that $u_\alpha\to u$ strongly (resp. weakly, weakly) and
$v_\alpha\to v$ weakly (resp. strongly, weakly).

\proof
The equivalence comes from the definition of the weak convergence and the fact that a net
$(u_\alpha)$ strongly converges to $u$ if and only if $\langle u_\alpha,v_\alpha\rangle_\alpha\to \langle u,v\rangle_\infty$ for every net $(v_\alpha)_\alpha$
weakly converging to $v\in L^2_\infty$. The "if" part is straightforward, for the "only if" we
see that for every net $(v_\alpha)$ strongly converging to $v$ we have $\langle u_\alpha,v_\alpha\rangle_\alpha\to \langle u,v\rangle_\infty$, which
implies the weak convergence of the net $(u_\alpha)$. Using now the hypothesis we get the convergence
of the net $\n{u_\alpha}_\alpha$ and thus the strong convergence of $(u_\alpha)$ by lemma \ref{wcvnet}.
\qed

\proclaim Proposition \theo. Let $(B_\alpha)$ be a strongly converging net to $B$ then
$$
\liminf_\alpha \n{B_\alpha}_{{\cal L}_\alpha}\geq \n{B}_{{\cal L}_\infty}
$$
and if the convergence is compact then it is an equality and 
$B$ is a compact operator as is its adjoint $B^*$.\label{ncvsnf}

\proof
Let $\epsilon>0$, there is $u\in {\cal L}^2_\infty$ such that $\n{u}_\infty=1$ and $\n{Bu}_\infty > \n{B}_{{\cal L}_\infty}-\epsilon$.
Take $(u_\alpha)_\alpha$ a net converging strongly to $u$. Then $\n{u_\alpha}_\alpha\to 1$, furthermore the strong convergence
of $(B_\alpha)$ implies that $\n{B_\alpha u_\alpha}_\alpha \to \n{Bu}_\infty$ thus
$$
\liminf_\alpha \n{B_\alpha}_{{\cal L}_\alpha} \geq \liminf_\alpha {\n{B_\alpha u_\alpha}_\alpha \over \n{u_\alpha}_\alpha } = \n{Bu}_\infty > \n{B}_{{\cal L}_\infty}-\epsilon {\rm .}
$$
Let suppose now that the convergence is compact. Take a net $(u_\alpha)$ such that $\n{u_\alpha}_\alpha=1$ and
$$
\lim_\alpha \bigl| \n{B_\alpha}_{{\cal L}_\alpha} - \n{B_\alpha u_\alpha}_\alpha \bigr| =0
$$
extracting a subnet if necessary we can suppose that the net $(u_\alpha)$ weakly converges to $u$.
By lemma \ref{wcvnet} we have $\n{u}_\infty\leq 1$, furthermore the compact convergence
 implies the strong convergence of $B_\alpha u_\alpha$ to $Bu$ thus
$$
\n{B}_{{\cal L}_\infty} \geq {\n{Bu}_\infty  \over \n{u}_\infty} \geq \n{Bu}_\infty = \lim_\alpha \n{B_\alpha u_\alpha}_\alpha = \lim_\alpha  \n{B_\alpha}_{{\cal L}_\alpha}
$$

Now let us prove that in the latest case, $B$ is compact. Let $(v_\beta)_{\beta\in{\cal B}}$ a net weakly converging
to $v$ in ${\cal L}^2_\infty$ then
$$
\langle u,Bv_\beta\rangle_\infty=\langle B^*u,v_\beta\rangle_\infty \to \langle B^*u,v\rangle_\infty = \langle u,Bv\rangle_\infty
$$
thus $Bv_\beta$ weakly converges to $Bv$. For every $\beta$ let $(u_{\alpha,\beta})$ be a strongly converging net
such that $\lim_\alpha u_{\alpha,\beta}=v_\beta$. For every $\beta$ the compact convergence of $(B_\alpha)$ implies the strong convergence
of $B_\alpha u_{\alpha,\beta}$ to $Bv_\beta$. Now let us take a net of positive numbers such that $\lim_\beta \epsilon(\beta)=0$,
then there is $\alpha(\beta)$ such that for every $\alpha\geq\alpha(\beta)$ we have
$$
\bigl| \n{B_\alpha u_{\alpha,\beta}}_\alpha-\n{Bv_\beta}_\infty \bigr| \leq \epsilon(\beta) {\rm .}
$$
Set $w_\beta=u_{\alpha(\beta),\beta}$ then $\lim_\beta w_\beta=v$ weakly and by the compact convergence we obtain the strong
convergence of $(B_{\alpha(\beta)}w_\beta)_\beta$ to $Bv$ but
$$
\lim_\beta \bigl| \n{B_{\alpha(\beta)}w_\beta}_{\alpha(\beta)} -\n{Bv_\beta}_\infty \bigr| =0
$$
which implies $\n{Bv_\beta}_\infty\to \n{Bv}_\infty$. We can conclude using lemma \ref{wcvnet}.
\qed
\subpartie{Convergence of spectral structures}
Here we see $L^2_\alpha$ as a Hilbert space. Then $A_\alpha$ and $A$ will be self-adjoint operators, $E_\alpha$ and
$E$ their respective spectral measure and $R_\mu$, $R$ their resolvents for $\mu$ in the resolvent
space. We want to study the links between the convergence of $(A_\alpha)$, $(E_\alpha)$ and $(R^\alpha_\mu)$. The following
theorem says that it is the same

\proclaim Theorem \theo. %
Let $(A_\alpha)$ and $A$ be self-adjoint operators $E_\alpha$, $E$ their spectral
measures and $R^\alpha_\mu$, $R_\mu$ their resolvents for $\mu$ in the resolvent space, then the following
assertions are equivalent
\item{1.}$R^\alpha_\mu \to R_\mu$ strongly (resp. compactly) for $\mu$ outside 
the union of the spectra of $A_\alpha$ and $A$.
\item{2.}$\varphi(A_\alpha)\to \varphi(A)$ strongly (resp. compactly) for every continuous function, with compact
support $\varphi:\R\to \C$
\item{3.}$\varphi_\alpha(A_\alpha)\to \varphi(A)$ strongly (resp. compactly) for every net $\{\varphi_\alpha:\R \to \C \}$ of continuous
functions vanishing at infinity and uniformly converging to $\varphi$ a continuous function vanishing
at infinity.
\item{4.}$E_\alpha\bigl(]\lambda,\mu]\bigr) \to E\bigl(]\lambda,\mu]\bigr)$ strongly (resp. compactly) for every pair of real
numbers outside the spectrum of $A$.
\item{5.}$\langle E_\alpha u_\alpha,v_\alpha\rangle_\alpha \to \langle Eu,v\rangle_\infty$ vaguely for every net of vectors $(u_\alpha)_{\alpha\in{\cal A}}$ and $(v_\alpha)_{\alpha\in{\cal A}}$
such that $u_\alpha\to u$ strongly (resp. weakly) and $v_\alpha\to v$ weakly.\label{sscve}

Let us recall that a quadratic form ${\cal Q}$ on a complex (resp. real) Hilbert space ${\cal H}$
comes from a sesquilinear (resp. bilinear) form, positive and symmetric 
${\cal E}: D({\cal E})\times D({\cal E}) \to \C$ (resp. $\R$) where $D(E)\in {\cal H}$ is a linear subspace
and ${\cal Q}(u)={\cal E}(u,u)$. Notice that ${\cal E}_1(u,v) = \langle u,v\rangle_{\cal H}+{\cal E}(u,v)$ for
every $u$ and $v\in D({\cal E})$ is also a sesquilinear (resp. bilinear), symmetric and positive form.
Thus $\bigl((D({\cal E}),{\cal E}_1\bigr)$ is a pre-Hilbert space. We say that ${\cal Q}$ is closed
if and only if $\bigl((D({\cal E}),{\cal E}_1\bigr)$ is a Hilbert space. In what follows, we will
not distinguish ${\cal Q}$ and the functional ${\cal E}$ defined by ${\cal E}(u)={\cal Q}(u)$ on
$D({\cal E})$ and ${\cal E}(u)=\infty$ on ${\cal H}\backslash D ({\cal E})$. In this context, ${\cal Q}$ is closed
if and only if ${\cal E}$ is lower semi-continuous as a function ${\cal E}:{\cal H}\to \overline{\R}$.

\proclaim Definition \theo. %
Let $({\cal E}_\alpha)$ be a net of closed quadratic forms, where
${\cal E}_\alpha$ is a closed quadratic form on $L^2_\alpha$ for every $\alpha\in {\cal A}$. We will say that
this net is {\rm asymptotically compact} if and only if for every net $(v_\alpha)_{\alpha\in{\cal A}}$ such that
$$
\limsup_\alpha {\cal E}_\alpha(v_\alpha) +\n{v_\alpha}_\alpha^2 < \infty
$$
there is a strongly converging sub-net.

Now a {\bf spectral structure} on a Hilbert space ${\cal H}$ over $\C$ (resp. $\R$) is
a family
$$
\Sigma=\{A,{\cal E},E,(T_t),(R_\zeta)\}
$$
where $A$ is a self-adjoint operator seen as the infinitesimal generator of the densely defined 
quadratic form ${\cal E}$ (such that $D({\cal E})=D(\sqrt A)$ and ${\cal E}(u,v)=\langle\sqrt Au,\sqrt Av\rangle_{\cal H}$
for every $u$ and $v$ in $D({\cal E}$) ), $E$ is its spectral measure, $(T_t)_{t\geq0}$ is a one parameter
semi-group of strongly continuous contractions ($T_t=e^{-tA}$, $t\geq0$) and $R_\zeta$ is a strongly continuous
resolvent ($R_\zeta=(\zeta-A)^{-1}$ for $\zeta\in \rho(A)$, where $\rho(A)$ is the resolvent set of $A$).
In what follows we will study a family of spectral structures $\Sigma_\alpha$ on $L^2_\alpha$, thus we will
have
$$
\Sigma_\alpha=\{A_\alpha,{\cal E}_\alpha,E_\alpha,(T_t^\alpha),(R_\zeta^\alpha)\}
$$

\proclaim Definition \theo. Let $(\Sigma_\alpha)_{\alpha\in{\cal A}}$ be a net with $\Sigma_\alpha$ a spectral structure on $L_\alpha^2$
and $\Sigma$ a spectral structure on $L^2_\infty$, we will say that the net $(\Sigma_\alpha)_\alpha$ strongly (resp. compactly)
converges to $\Sigma$ if and only if one of the conditions of theorem \ref{sscve} is satisfied.

\proclaim Proposition \theo. Let $(\Sigma_\alpha)_{\alpha\in{\cal A}}$ of spectral structures strongly
converging to $\Sigma$ then for any net $(v_\alpha)_\alpha$ weakly converging to $v$ we have
$$
{\cal E}(v)\leq \liminf_\alpha {\cal E}_\alpha(v_\alpha)
$$
Furthermore, if the net $(\Sigma_\alpha)_{\alpha\in{\cal A}}$ converges compactly, then the net
of quadratic forms $({\cal E}_\alpha)_\alpha$ is asymptotically compact.\label{thdefstrspe}

\proof
Assume that the net of resolvents $(R_\lambda^\alpha)$ is strongly convergent and write
$$
a^\lambda_\alpha(u,v)=-\lambda\scal{u-\lambda R_\lambda^\alpha u,v}_\alpha
$$ 
(the Deny-Yosida  approximation of
bilinear form associated to ${\cal E}_\alpha$), then the net $(a^\lambda_\alpha(u,u))$ converges to ${\cal E}_\alpha(u)$
increasing when $\lambda\to -\infty$ (see Mosco  \cite{mosco} 1.(i)).
From the assumption it easy to see that for $(u_\alpha)$ and $(v_\alpha)$ converging
strongly to $u$ and weakly to $v$ respectively
$$
\lim_\alpha a^\lambda_\alpha(u_\alpha,v_\alpha)=-\lambda\scal{u-\lambda R_\lambda u,v}_\infty=a^\lambda(u,v)
$$ 
we recall that (see Dal Maso \cite{Dmaso} proposition 12.12)
$$
a^{\lambda}(u,u)\geq a^\lambda(v,v)+2\lambda\scal{v-\lambda R_\lambda v,u-v}_\infty 
$$
hence for any net $v_\alpha$ weakly converging to $u$ and $w_\alpha$ a strongly converging net to $u$
we have
$$
{\cal E}_\alpha(v_\alpha)\geq a_\alpha^\lambda(v_\alpha,v_\alpha)\geq a^\lambda_\alpha(w_\alpha,w_\alpha)+2\lambda\scal{w_\alpha-\lambda R_\lambda^\alpha w_\alpha,v_\alpha-w_\alpha}
$$
thus $\liminf_\alpha {\cal E}_\alpha(v_\alpha) \geq a^\lambda(u,u)$ for any $\lambda<0$, now taking $\lambda\to -\infty$ we can conclude that
$\liminf_\alpha {\cal E}_\alpha(v_\alpha) \geq {\cal E}(u)$.

  Now assume that $(\Sigma_\alpha)$ compactly converges and let $(u_\alpha)_{\alpha\in{\cal A}}$ be a net such that
$$
\sup_\alpha \bigl( {\cal E}_\alpha(u_\alpha) +\n{u_\alpha}_\alpha^2 \bigr) \leq M < \infty {\rm\ .}
$$
Taking a sub-net if necessary we can suppose that $(u_\alpha)_\alpha$ weakly converges to $u$.
Let $\rho>0$ be out of $A_\infty$'s spectrum. As
$$
\int_{]\rho,\infty[} d\scal{E_\alpha u_\alpha,u_\alpha}_\alpha \leq \frac{1}{\rho}\int_{]\rho,\infty[} \lambda d\scal{E_\alpha(\lambda)u_\alpha,u_\alpha}_\alpha \leq \frac{{\cal E}_\alpha(u_\alpha)}{\rho}\leq \frac{M}{\rho}
$$
we have
$$
\n{u_\alpha}_\alpha^2 \leq \int_{[0,\rho]} d\scal{E_\alpha u_\alpha,u_\alpha}_\alpha + \frac{M}{\rho}
$$
and the compact convergence implies $\lim_\alpha \int_{[0,\rho]} d\scal{E_\alpha u_\alpha,u_\alpha}_\alpha = \int_{[0,\rho]} d\scal{Eu,u}_\infty$ hence
$$
\limsup_\alpha \n{u_\alpha}_\alpha^2 \leq \int_{[0,\rho]} d\scal{Eu,u}_\infty + \frac{M}{\rho} \leq \n{u}_\infty^2 + \frac{M}{\rho}
$$
and taking $\rho\to \infty$ we get
$$
\limsup_\alpha \n{u_\alpha}_\alpha^2 \leq \n{u}_\infty^2
$$
finally we deduce the strongly convergence of the net $(u_\alpha)$ using lemma \ref{wcvnet}.
\qed

The main reason we introduced all these convergences is the following
theorem, the proof of which we postpone to avoid drowning the reader in
too many technical details.
\proclaim Theorem \theo. Let $\Sigma_\alpha\to \Sigma$ compactly and suppose that all
resolvents $R^\alpha_\zeta$ are compact. Let $\lambda_k$ (resp. $\lambda_k^\alpha$) be the $k^{\rm th}$
eigenvalue of $A$ (resp. $A_\alpha$) with multiplicity. We take $\lambda_k=+\infty$ if
$k>\dim L^2_\infty +1$ when $\dim L^2_\infty< \infty$ and $\lambda_k^\alpha=+\infty$ if
$k>\dim L^2_\alpha +1$ when $\dim L^2_\alpha< \infty$. Then for every $k$
$$
\lim_\alpha \lambda_k^\alpha=\lambda_k
$$
Furthermore let $\{\varphi_k^\alpha\mid k=1,\ldots,\dim L_\alpha^2\}$ be an orthonormal bases of $L^2_\alpha$ such that
$\varphi_k^\alpha$ is an eigenfunction of $A_\alpha$ for $\lambda_k^\alpha$. Then there is a sub-net such that
for all $k\leq \dim L^2_\infty$ the net $(\varphi_k^\alpha)_\alpha$ strongly converges to the eigenfunction 
$\varphi_k$ of $A$ for the eigenvalue $\lambda_k$, and such that the family
 $\{\varphi_k\mid k=1,\ldots,\dim L_\alpha^2\}$ is an orthonormal basis of $L^2_\infty$.\label{cvssets}

\partie{Proof of Theorem 1}\label{proof1}
\subpartie{Homogenisation of the Laplacian }\label{homogelap}
In this section we are going to built the operator $\Delta_\infty$ of theorem \ref{theo1}.
we remind the reader that $D_f$ is a fundamental domain, we then begin by
taking $\chi^i$ as the unique periodic solution (up to an additive constant)
of
$$
\Delta\chi^i = \Delta x_i \ {\rm on}\ D_f
$$
The operator $\Delta_\infty$ is then defined by
$$
\Delta_\infty f = -{1 \over {\rm Vol}(g)} \Biggl(\int_{D_f} g^{ij}-g^{ik} {\partial \chi^j \over \partial y_k}~d\mu_g \Biggr)%
{\partial^2f\over \partial x_i  \partial x_j} \numeq
$$
Now let us write $\eta_j(x)=\chi^j(x)-x_j$ the induced harmonic function and
$$
q^{ij}= {1 \over {\rm Vol}(g)} \Biggl(\int_{D_f} g^{ij}-g^{ik} {\partial \chi^j \over \partial y_k}~d\mu_g \Biggr)
$$
we can notice that the $d\eta_i$ are harmonic $1$-forms on the torus.
It is not difficult now to show that

\proclaim Proposition \theo.
Let $\langle\cdot,\cdot\rangle_2$ be the scalar product induced on $1$-forms by the Riemannian
metric $g$. Then
$$
q^{ij}= {1\over {\rm Vol}(g) } \langle d\eta_i,d\eta_j\rangle_2 =q^{ji}
$$
thus $\Delta_\infty$ is an elliptic operator.

In fact we can say more, $(q^{ij})$ induces a scalar product on harmonic $1$-forms
(whose norm will be written $\n{\cdot}_2$) and then to $H^1(\T,\R)$. Indeed,
as mentioned earlier, we can see the $(d\eta_i)$ as $1$-forms over the torus.
Being a free family they can be seen as a basis of $H^1(\T,\R)$ (Hodge's theorem).
Thus by duality this yields also a scalar product $(q_{ij})$ over $H_1(\T,\R)$ 
(whose induced norm will be written $\n{\cdot}^*_2$). The question naturally arising
is to know the link between this norm and the stable norm. To see this we have to
go back on $H^1(\T,\R)$. Indeed the stable norm is the dual of the norm
obtained by quotient of the sup norm on $1$-forms (see Pansu \cite{pansu2} lemma 17),
which we write $\n{\cdot}_\infty^*$, and the norm $\n{\cdot}_2$ comes from the
normalised $L^2$ norm. Thus mixing the H{\"o}lder inequality and the Hodge-de~Rham
theorem we get :
\proclaim Proposition \theo. For every $1$-form $\alpha$ we have
$$
\n{\alpha}_2 \leq \n{\alpha}_\infty^*
$$
thus by duality, for every $\gamma\in H_1(\T,\R)$ we have
$$
\n{\gamma}_\infty \leq \n{\gamma}^*_2
$$
in other words the unit ball of $\n{\cdot}^*_2$ is included in $B_\infty(1)$.

To finish this section, let us remark that the manifold $H_1(\T,\R)/ H_1(\T,\Z)$
with the flat metric induced by $\n{\cdot}^*_2$ is usually called the Jacobi manifold
or the Albanese torus of $(\T,g)$.
\subpartie{Asymptotic compactness}
Let us now define the various functional spaces involved. For $\rho\in \overline{\R}$,
$L^2\bigl(B_\rho(1),d\mu_\rho\bigr)$ will be the space of square integrable functions
over the ball $B_\rho(1)$, which is a Hilbert space with the scalar product
$$
(u,v)_\rho=\int_{B_\rho(1)} uv~d\mu_\rho
$$
whose norm will be $|\cdot|_\rho$. $H^1_{\rho,0}\bigl(B_\rho(1)\bigr)$ will be the closure
of $C^\infty\bigl(B_\rho(1)\bigr)$ functions with compact support, in $H^1_\rho\bigl(B_\rho(1)\bigr)$
for the norm $\n{\cdot}_\rho$ defined by
$$
\n{v}_\rho^2= |v|_\rho^2 + \sum_{i=1}^n \Biggv{{\partial v\over \partial x_i}}_\rho^2
$$
and with
$$
H^1_\rho\bigl(B_\rho(1)\bigr) = \biggl\{ v \biggm|v,{\partial v \over \partial x_1},\ldots,{\partial v \over \partial x_n} \in L^2\bigl(B_\rho(1),d\mu_\rho\bigr) \biggr\}
$$

For all that follows, $V_\rho$ will be a closed sub-space such that
$$
H^1_{\rho,0}\bigl(B_\rho(1)\bigr)\subset V_\rho\subset H_\rho^1\bigl(B_\rho(1)\bigr)
$$
Thus we can define a spectral structure on $L^2_\rho$ by expanding the Laplacian
defined on $V_\rho$ on $L^2_\rho$. If $V_\rho=H^1_{\rho,0}\bigl(B_\rho(1)\bigr)$ we deal with
the Dirichlet problem, and if $V_\rho=H^1_\rho\bigl(B_\rho(1)\bigr)$ we then deal with
the Neumann problem.
We then put the following norm on $V_\rho$ :
$$
\n{v}_{\rho,0}^2 = |v|_\rho^2 + (v,\Delta_\rho v)_\rho 
$$
we then have

\proclaim Lemma \theo.
Let $(u_\rho)_\rho$ be a net with $u_\rho \in V_\rho$ for every $\rho$, if there is
a constant $C$ such that for all $\rho>0$ we have
$$
\n{u_\rho}_{\rho,0} \leq C
$$
then there is a strongly converging sub-net in ${\cal L}^2$.\label{lemapivot}

\proof
Let $B=\cup_\rho B_\rho(1)$ we are going to show
that the strong convergence in $L^2(B,\mu_\infty)$ implies the strong convergence
in ${\cal L}^2$. Then the compact embedding of $H^1_\infty\bigl(B\bigr)$ in $L^2\bigl(B,\mu_\infty\bigr)$
will conclude the proof.

Let us first notice that the periodicity gives the existence of two constant
$\alpha$ and $\beta$ such that
$$
\alpha |v|_\infty \leq |v|_\rho \leq \beta |v|_\infty {\rm\ .}
$$
Let us start by taking a net $(u_\rho)$ strongly converging in $L^2(B,\mu_\infty)$ to $u_\infty$
we also assume $u_\rho \in V_\rho$ for every $\rho$, because it is all we need. Now let
$c_p \in C_0^\infty\bigl( B_\infty(1) \bigr)$ be a sequence of functions strongly converging to $u_\infty$
and take $p$ large enough for the support of $u_p$ to be in $B_\rho(1)$. We have
$$
|c_p-u_\rho|_\rho \leq \beta |c_p-u_\infty|_\infty +\beta |u_\infty-u_\rho|_\infty
$$
now let  $\varepsilon>0$ then for $p$ large enough $\beta |c_p-u_\infty|_\infty \leq \varepsilon$. We fix $p$
large enough and take $\rho$ large enough for the second term to converge to 0.

Now to conclude observe that from the assumptions the net $(u_\rho)$ is bounded in
$H^1_\infty\bigl(B\bigr)$, hence using the compact embedding of $H^1_\infty\bigl(B\bigr)$
in $L^2(B,\mu_\infty)$ we can extract a strongly converging net in $L^2(B,\mu_\infty)$ and by what
we just did in ${\cal L}^2$.
\qed

\subpartie{Compact convergence of the resolvents}
Let $\lambda>0$ and $G_\lambda^\rho$ be the operator from $L^2_\rho$ to $V_\rho\subset L^2_\rho$
such that
$$
a^\rho_\lambda(G_\lambda^\rho f,\phi) = (f,\phi)_\rho\quad \forall \phi \in V_\rho {\rm.} \numeq \label{predemo1}
$$
where 
$$
a^\rho_\lambda(u,v)= \int_{B_\rho(1)} g_\rho^{ij}~\partial_iu\cdot~\partial_jv~d\mu_\rho + \lambda(u,v)_\rho
$$
We want to show that the net of operators $(G_\lambda^\rho)$ converges compactly to
$G_\lambda$ the operator corresponding to the homogenised problem:
$$
a^\infty_\lambda(G_\lambda f,\phi)=(f,\phi)_\infty  \quad \forall \phi \in V_\infty  \numeq \label{predemo2}
$$
with $(f,\phi)_\infty=\int_{B_\infty(1)} f\phi~d\mu_\infty$ and
$$
a^\infty_\lambda(u,v)=\int_{B_\infty(1)} q^{ij}~\partial_iu~\partial_jv~d\mu_\infty+\lambda(u,v)_\infty
$$
in other word we want to show the following theorem

\proclaim Theorem \theo.
  For every $\lambda<0$, the net of resolvents $(R_\lambda^\rho)_\rho$ of the Laplacian $(\Delta_\rho)$ converges 
compactly to $R_\lambda^\infty$, the resolvent of $\Delta_\infty$ from the homogenised problem.
Thus the net $(\Sigma_\rho)$ compactly converges to $\Sigma_\infty$.\label{cvres}

\proof
This comes from the fact that $R_\lambda^\rho=-G_{-\lambda}^\rho$ and $R_\lambda^\infty=-G_{-\lambda}$.

\noindent{\bf First step :}

Let $f_\rho$ be a weakly convergent net to $f$ in
${\cal L}^2$, thus from \ref{wcvnet} this net is uniformly bounded in ${\cal L}^2$ and in
$V_\rho'$, the dual space of $V_\rho$.

Let $f_\rho \in V_\rho$ then by \ref{predemo1} we have :
$$
\alpha \n{G^\rho_\lambda f_\rho}^2_{\rho,0} \leq (f_\rho,G_\lambda^\rho f_\rho)_\rho \leq %
K\n{f_\rho}_{V_\rho'}\n{G_\lambda^\rho f_\rho}_{\rho,0}
$$
thus
$$
\n{G_\lambda^\rho f_\rho}_{\rho,0} \leq C \n{f_\rho}_{V_\rho'}
$$
the net $(G_\lambda^\rho f_\rho)$ being uniformly bounded for the norms $\n{\cdot}_{\rho,0}$, using lemma \ref{subwcvnet}
 there is a subnet strongly converging in ${\cal L}^2$.  i.e.
$$
u_\rho=G_\lambda^\rho f_\rho \to u^*_\lambda {\rm\ strongly\ in\ } {\cal L}^2 \numeq \label{demo1}
$$
Furthermore  $P_\rho=(g^{ij}_\rho)\nabla G_\lambda^\rho f_\rho$ is also bounded in
${\cal L}^2$ thus there is a subnet of the net $P_\rho$ 
weakly converging in ${\cal L}^2$ to $P^*_\lambda\in L^2_\infty$.  For any
$\phi_\infty\in L^2_\infty $ let $\phi_\rho$ be a strongly converging net to
$\phi_\infty$ in ${\cal L}^2$ then
$$
  \eqalign{%
  \int_{B_\rho(1)} P_\rho\!\cdot\!\nabla\phi_\rho~d\mu_\rho +\lambda(G_\lambda^\rho f_\rho,\phi_\rho)_\rho&= (f_\rho,\phi_\rho)_\rho \to \cr %
  \int_{B_\infty(1)} P^*_\lambda\!\cdot\!\nabla\phi_\infty~d\mu_\infty+\lambda(u^*_\lambda,\phi_\infty)_\infty&=(f,\phi_\infty)_\infty{\rm .}\cr}\numeq \label{demo1.1}
$$
Thus it is enough to show that $P^*_\lambda=\bigl(q^{ij}\bigr)\nabla u^*_\lambda$ on
$B_\infty(1)$ because it induces $u^*_\lambda=G_\lambda f$.

\noindent{\bf Second step : } 

We first take $\chi^k(y)$ (see \ref{homogelap}) such that ${\cal M}(\chi^k)= 0$ and we define
$$
w_\rho(x)=x_k- \frac{1}{\rho} \chi^k(\rho x) \numeq \label{demo2}
$$
for every $k=1,\dots,d_1$.
Then
$$
w_\rho \to x_k \ {\rm  strongly\ in }\  {\cal L}^2.\numeq \label{demo2.1}
$$  
an by construction  of $\chi^k$ (see \ref{homogelap})  we have
$$
 -\partial_i\bigl(F(g_\rho,X) g^{ij}_\rho~\partial_jw_\rho \bigr) =0 \ {\rm on }\ B_\rho(1) {\rm .} \numeq \label{demo3}
$$
We multiply this equation by a test function 
$\phi \in V_\rho$ and after an integration we get 
$$
\int_{B_\rho(1)}g^{ij}_\rho~\partial_jw_\rho~\partial_i\phi~d\mu_\rho=0  \numeq \label{demo4}
$$ 
Let $\varphi \in C^\infty_0(B_\infty(1))$ (notice that for $\rho$ large enough 
the support of $\varphi$ will be in $B_\rho(1)$ ) and $\phi=\varphi w_\rho$ which we
put into the equation \ref{predemo1} and into the equation \ref{demo4} we put $\phi=\varphi u_\rho$ 
(see \ref{demo1}), and then we subtract the results
$$\displaylines{%
\qquad \int_{B_\rho(1)} g^{ij}_\rho\bigl( {\partial_ju_\rho}~{\partial_i\varphi}~w_\rho - \partial_jw_\rho~\partial_i\varphi~u_\rho \bigr)~d\mu_\rho \hfill \cr
\hfill =\int_{B_\rho(1)} f_\rho w_\rho\varphi~d\mu_\rho -\lambda \int_{B_\rho(1)}\varphi u_\rho w_\rho~d\mu_\rho \qquad \nume \label{demo5} }
$$
Now let  $\rho \to \infty $ in \ref{demo5}, all
terms converge because they are product of one strongly converging net and one weakly 
converging net in ${\cal L}^2$.
More precisely
    \item{---} $P_\rho$ defined $P_{\rho,i}=g^{ij}_\rho \partial_ju_\rho$ weakly converges 
  to $P^*_\lambda$ in ${\cal L}^2$ following \ref{demo1.1};
    \item{---} $\partial_i\varphi w_\rho$ strongly converges to
  $\partial_i\varphi x_k$ in ${\cal L}^2$ from \ref{demo2.1}
    \item{---} $g^{ij}_\rho \partial_iw_\rho$ is $D_f/\rho$-periodic and weakly converges in ${\cal L}^2$ towards its mean
value
  $$
  q^{jk}={\cal M}\biggl(g^{ij}(y)\Bigl(\delta_{ik}-\partial_i \chi^k(y)\Bigr)\biggr)
  $$
    \item{---} $\partial_j\varphi u_\rho$ strongly  converges to
  $\partial_j\varphi u^*_\lambda$ by \ref{demo1}, because $\varphi$ has compact support.
    \item{---} Now for the right side, $w_\rho$ strongly converges as $u_\rho$ does 
and  $f_\rho$ weakly converges to $f$.

To summarise \ref{demo5} converges to (we write $P^*_{\lambda,i}$ the coordinates of $P_\lambda^*$)
$$
\int_{B_\infty(1)} \bigl(P^*_{\lambda,j}x_k - q^{jk} u^*_\lambda\bigr)\partial_j\varphi~d\mu_\infty= \int_{B_\infty(1)} fx_k \varphi~d\mu_\infty -\lambda \int_{B_\infty(1)}\varphi u^*_\lambda x_k~d\mu_\infty \numeq 
\label{demo6}
$$
furthermore if we put into equation \ref{demo1.1}, $\phi_\infty=\varphi x_k$ it gives
$$
 \int_{B_\infty(1)}fx_k\varphi~d\mu_\infty-\lambda\int_{B_\infty(1)}\varphi u^*_\lambda x_k~d\mu_\infty = \int_{B_\infty(1)} P_{\lambda,j}^* \partial_j(\varphi x_k)~d\mu_\infty  \numeq \label{equa}
 $$
and by mixing \ref{demo6} and \ref{equa} we get
for every $\varphi \in {\cal C}^\infty_c(B_\infty(1))$ the following equality : 
$$
\int_{B_\infty(1)}  \bigl(P^*_{\lambda,j}x_k - q^{jk} u^*\bigr)\partial_j\varphi~d\mu_\infty =\int_{B_\infty(1)}P_{\lambda,j}^* \partial_j(\varphi x_k)~d\mu_\infty
$$
which in terms of distribution can be translated into :
$$
-\sum_{j=1}^{d_1}\partial_j\bigl(P_{\lambda,j}^*x_k - q^{jk} u^*_\lambda\bigr) =-\sum_{j=1}^{d_1} \partial_jP_j^*~x_k
\iff  P_{\lambda,k}^*=\sum_{j=1}^{d_1}q^{jk}~\partial_ju^*_\lambda
$$
which allow us to conclude that $u^*_\lambda=G_\lambda f$. \qed

It is now easy to finish the proof of theorem \ref{theo1}, it comes from theorem \ref{cvres}
and theorem \ref{cvssets}.
\partie{$\sectionfont \Gamma$-convergence of quadratic forms}\label{topo2}
\subpartie{$\subsectionfont \Gamma$  and  Mosco-convergence of quadratic forms}
We are now going to give a definition of $\Gamma$-convergence adapted to our problem.

\proclaim Definition {\theo.} ($\Gamma$-convergence).
We say that a net $\{ F_\alpha : L^2_\alpha \to \overline{\R}\}_{\alpha\in {\cal A}}$ of functions
{\rm $\Gamma$-converges to $F : L^2_\infty\to \overline{\R}$} if and only if the
following assertions are satisfied :
    \item{(F1)} For any net $(u_\alpha)_{\alpha\in {\cal A}}\in L^2_\alpha$
  strongly converging to $u\in L_\infty^2$ in ${\cal L}^2$ we have
  $$
  F(u) \leq \liminf_\alpha F_\alpha(u_\alpha);
  $$
    \item{(F2)} For every $u\in L_\infty^2$ there is a net $(u_\alpha)_{\alpha\in
    {\cal A}}\in L^2_\alpha$ strongly converging to  $u$ in
  ${\cal L}^2$ such that
  $$
  F(u) = \lim_\alpha F_\alpha (u_\alpha).\label{gammaconv}
  $$

\noindent{\bf Remark.}
This is slightly different from the definition \ref{gamacv1}, which is the usual one.
By taking $F_\alpha$ infinite outside of $L_\alpha^2$ in ${\cal L}^2$ we get back (in some way) the usual
definition (see the introduction of \cite{Dmaso}).

Let us summarise some properties satisfied by this convergence.

\proclaim Lemma \theo.
      \item{(a)} Let $\{F_\alpha :L^2_\alpha \to \overline{\R}\}_{\alpha\in {\cal A}}$ be a net
       of functions $\Gamma$-converging to a function $F:L_\infty^2\to\overline{\R}$, 
        then $F$ is lower semi-continuous.
      \item{(b)} Let $({\cal E}_\alpha)_{\alpha\in {\cal A}}$ be a net of
    quadratic forms ${\cal E}_\alpha$ on $L^2_\alpha$ $\Gamma$-converging to a
    function $F:L_\infty^2\to \overline{\R}$, then $F$ can be identified with a
    quadratic form on $L_\infty^2$.\label{fquadbis}

There is also the following result, concerning compactness

\proclaim Theorem \theo.
  From every net $({\cal E}_\alpha)_{\alpha\in {\cal A}}$ of quadratic forms 
  ${\cal E}_\alpha$ on $L^2_\alpha$ we can extract a  $\Gamma$-converging subnet,
  whose limit is a quadratic form on $L_\infty^2$. 

\noindent{\bf Remark.}
  This theorem is true for a wider variety of functions, with
some restrictions on  $\{ L^2_\alpha \}_{\nu \in {\cal A}}$. Of course the limit in
that case is not always a quadratic form. Here it is lemma \ref{fquadbis}
which gives information on the limit.

\proclaim Definition {\theo.} (Mosco topology).
We say that a net $({\cal E}_\alpha)_{\alpha\in {\cal A}}$ of quadratic forms ${\cal E}_\alpha$ on $L^2_\alpha$ 
{\rm Mosco-converges} to the quadratic form ${\cal E}$ on $L_\infty^2$ if condition (F2) 
of  definition \ref{gammaconv} and (F1') are satisfied:
    \item{(F1')} For any $(u_\alpha)_{\alpha\in {\cal A}},\ u_\alpha\in L^2_\alpha$ 
 weakly converging net to $u\in L_\infty^2$ in ${\cal L}^2$ we have
  $$
  {\cal E}(u) \leq \liminf_\alpha {\cal E}_\alpha(u_\alpha)
  $$
The induced topology is called the  {\rm Mosco  topology}.\label{moscocv}

It is obvious that the Mosco-convergence induces the $\Gamma$-convergence, thus this
topology is stronger. Let us now define one last convergence :

\proclaim Definition {\theo.} (Compact $\Gamma$-convergence).
We say that a net $({\cal E}_\alpha)_{\alpha\in {\cal A}}$ {\rm $\Gamma$-converges compactly} to ${\cal E}$
  if ${\cal E}_a \to {\cal E}$ in the Mosco topology an if $({\cal E}_a)_{a\in{\cal A}}$ 
  is asymptotically compact.

Let us precise how the Mosco and the $\Gamma$ topologies are linked :
\proclaim Lemma \theo.
  Let us suppose $({\cal E}_\alpha)_{\alpha\in {\cal A}}$ asymptotically
  compact then $({\cal E}_\alpha)_{\alpha\in {\cal A}}$ $\Gamma$-converges to
  ${\cal E}$ is and only if $({\cal E}_\alpha)_{\alpha\in {\cal A}}$ Mosco-converges 
  to ${\cal E}$.\label{maequg}

\proof
We just need to show that the $\Gamma$-convergence implies the condition
$(F1')$ from definition \ref{moscocv}. We proceed {\sl ad absurdum} and suppose 
that there is a weakly converging net $(u_\alpha)$ such that $\liminf_\alpha {\cal E}_\alpha(u_\alpha) < {\cal E}(u)$.
Taking a subnet if necessarily we can suppose $\lim {\cal E}_\alpha(u_\alpha) < {\cal E}(u)$ 
thus we also have $ \limsup_\alpha {\cal E}_\alpha(u_\alpha) + \n{u_\alpha}_\alpha^2 < +\infty $. 
The asymptotic compactness is obviously inherited by a subnet thus we can extract
a strongly converging  subnet $u_{\alpha(\beta)}$.
The $\Gamma$-convergence being also inherited by a subnet of  ${\cal E}_\alpha$ we finally get
$$
\lim {\cal E}_\alpha(u_\alpha)=\lim_\beta{\cal E}_{\alpha(\beta)}(u_{\alpha(\beta)}) \geq {\cal E}(u)
$$
which is absurd.
\qed

\subpartie{$\subsectionfont \Gamma$-convergence and spectral structures}

The following theorem explains how are related the convergence of spectral
structures and the Mosco-convergence.
\proclaim Theorem \theo.
  Let $(\Sigma_\alpha)$ be a net of spectral structures on $(L^2_\alpha)$ and $\Sigma$ a spectral structure
on $L^2_\infty$ then  $\Sigma_\alpha \tend \Sigma$ strongly (resp. compactly) if and only if
${\cal E}_\alpha$ Mosco-converges (resp. $\Gamma$-converges compactly) to ${\cal E}$.

\proof
We are going to prove the equivalence between the strong (resp. compact) convergence of
resolvents and the Mosco-convergence (resp. $\Gamma$-convergence compact) of the energies.

Let us begin by assuming the Mosco-convergence of the net  $({\cal E}_\alpha)$.
We need to show that for every $z\in L^2_\infty$ and any net $(z_\alpha)$ strongly converging to $z$
the net $u_\alpha=-R^\alpha_\lambda z_\alpha$ strongly converges to $u=-R_\lambda z$.
 First let us notice that the vector $u$ is the unique minimiser of
$$
v \mapsto {\cal E}(v)-\lambda\n{v}_\infty^2-2\scal{z,v}_\infty
$$
we can characterise the same way  $u_\alpha$ for every $\alpha$.

As an operator of $L^2_\alpha$, $R_\lambda^\alpha$ is bounded by ${-\lambda}^{-1}$.
Thus the net $(u_\alpha)$ is bounded and we can extract a weakly converging sub-net, still
written $(u_\alpha)$, with limit $\tilde{u}$. Now from condition 
(F2) for every $v\in L^2_\infty$ there is a net strongly converging to it such that
$\lim_\alpha {\cal E}_\alpha(v_\alpha)={\cal E}(v)$. But for every $\alpha$
$$
{\cal E}_\alpha(u_\alpha)-\lambda\n{u_\alpha}^2_\alpha-2\scal{z_\alpha,u_\alpha}_\alpha\leq{\cal E}_\alpha(v_\alpha)-\lambda\n{v_\alpha}^2_\alpha-2\scal{z_\alpha,v_\alpha}_\alpha \numeq \label{mequres1}
$$
thus taking the limit in $\alpha\in{\cal A}$ we get thanks to condition (F1') of definition\ref{moscocv}
and the fact that for any weakly convergent net $\n{\tilde{u}}_\infty\leq \liminf_\alpha \n{u_\alpha}_\alpha$ 
(remember that $\lambda<0$)
$$
{\cal E}(\tilde{u})-\lambda\n{\tilde{u}}_\infty^2-2\scal{z,\tilde{u}}_\infty \leq{\cal E}(v)-\lambda\n{v}_\infty^2-2\scal{z,v}_\infty
$$
which implies $\tilde{u}=-R_\lambda z$. Due to $u$'s unicity,
we conclude that $(u_\alpha)$ weakly converges to $u$. 
Let us prove that $\n{u_\alpha}_\alpha$ converges to $\n{u}_\infty$.
In that aim take a strongly convergent net $v_\alpha$ to $v$ such that $\lim_\alpha {\cal E}_\alpha(v_\alpha)={\cal E}(u)$, 
and take a new look at inequality \ref{mequres1} :
$$
{\cal E}_\alpha(u_\alpha)- \lambda\n{u_\alpha+z_\alpha/\lambda}^2_\alpha \leq {\cal E}_\alpha(v_\alpha)- \lambda\n{v_\alpha+z_\alpha/\lambda}^2_\alpha
$$
using (F1') once again we find
$$
 {\cal E}(v)-\lambda\limsup_\alpha \n{u_\alpha+z_\alpha/\lambda}^2_\alpha \leq {\cal E}(v) - \lambda\n{u+z/\lambda}_\infty^2
$$
thus $\n{u_\alpha+z_\alpha/\lambda}_\alpha^2\to\n{u+z/\lambda}_\infty^2$ which implies the strong convergence of $(u_\alpha+z_\alpha/\lambda)_\alpha$ and
the strong convergence of $(z_\alpha)$ induces  the strong convergence of $(u_\alpha)$.

We shall now study the compact $\Gamma$-convergence.
Let us take a weakly convergent net $w_\alpha$ to $w$ and let $u_\alpha=-R^\alpha_\lambda w_\alpha$, then the net $u_\alpha$
is still bounded.
Swapping $z_\alpha$ with $w_\alpha$ in \ref{mequres1} we get that $\limsup_\alpha{\cal E}_\alpha(u_\alpha)$ is bounded,
and thanks to  the asymptotic compactness we can extract a strongly convergent sub-net
with $\tilde{u}$ its limit. Putting this in \ref{mequres1}, with $z_\alpha=v_\alpha$ where $(v_\alpha)$ a
strongly converging net to $v$ we get
$$
{\cal E}(\tilde{u})-\lambda\n{\tilde{u}}^2-2\scal{w,\tilde{u}} \leq{\cal E}(v)-\lambda\n{v}^2-2\scal{w,v}
$$
thus $\tilde{u}=-R_\lambda w$. 
Once again, thanks to unicity, we conclude that $R_\lambda^\alpha w_\alpha$ strongly converges to $R_\lambda w$.

{\bf Reciprocally} assume that for every $\lambda<0$ the net $R_\lambda^\alpha$ strongly converges to $R_\lambda$.
In what follows $(u_\alpha)$ will be a strong convergent net to $u$.

\noindent{\bf Condition (F1') :} already done, see proposition \ref{thdefstrspe}.

\noindent{\bf Condition (F2) :} Extract a sub-net $\lambda_\alpha\to -\infty$ such that
$$
{\cal E}(u,u) \geq \lim_\lambda \lim_\alpha a_\alpha^\lambda(u_\alpha,u_\alpha) \geq \lim_\alpha a_\alpha^{\lambda_\alpha}(u_\alpha,u_\alpha)
$$
take $w_\alpha=\lambda_\alpha R_{\lambda_\alpha}^\alpha u_\alpha$ for every $\alpha$ and notice that
$$\eqalign{
a^\lambda_\alpha(u_\alpha,u_\alpha)=&-\lambda\scal{u_\alpha-\lambda R_\lambda^\alpha u_\alpha,u_\alpha}_\alpha-\lambda\scal{u_\alpha-\lambda R_\lambda^\alpha u_\alpha,-\lambda R^\alpha_\lambda u_\alpha}_\alpha\cr
            & +\lambda\scal{u_\alpha-\lambda R_\lambda^\alpha u_\alpha,-\lambda R^\alpha_\lambda u_\alpha}_\alpha \cr
=&-\lambda\n{u_\alpha-\lambda R_\lambda^\alpha u_\alpha}^2+\lambda^2\scal{u_\alpha-\lambda R_\lambda^\alpha u_\alpha,-R^\alpha_\lambda u_\alpha}_\alpha \cr
=&-\lambda\n{u_\alpha-\lambda R_\lambda^\alpha u_\alpha}^2+{\cal E}_\alpha(\lambda R^\alpha_\lambda u_\alpha)\cr}
$$
indeed if  $a_\alpha$ is the bilinear form corresponding to ${\cal E}_\alpha$ then
$R_\lambda^\alpha u_\alpha$ can be seen as the sole element such that 
$$
a_\alpha(-R_\lambda^\alpha u_\alpha,v_\alpha)-\lambda\scal{-R_\lambda^\alpha u_\alpha,v_\alpha}_\alpha=\scal{u_\alpha,v_\alpha}_\alpha, \quad \forall v_\alpha \in D({\cal E}_\alpha)
$$
hence 
$$
a_\alpha^{\lambda_\alpha}(u_\alpha,u_\alpha)={\cal E}_\alpha(w_\alpha,w_\alpha)-\lambda_\alpha\n{u_\alpha-w_\alpha}^2
$$
which implies $w_\alpha \to u$ strongly in ${\cal L}^2$ and
$$
{\cal E}(u,u)\geq\limsup_{\alpha \to +\infty} {\cal E}_\alpha(w_\alpha,w_\alpha)
$$

For the compact convergence case it suffices to prove the asymptotic compactness,
but it has already been done in the proof of proposition \ref{thdefstrspe}.
\qed
\partie{Proof of Theorem 2}\label{proof2}
      
The convergence of the eigenvalue is given by theorem \ref{theo1}.
Hence it remains to bound the asymptotic $\lambda_1$ (i.e. the limit) and
characterise the equality.
The proof we propose consists in finding an upper bound of $\lambda_1\bigl(B_g(\rho)\bigr)$ 
for every $\rho$ using a function depending of the distance from the centre of
the ball. We then use the simple convergence of the distances $(d_\rho)$ to the
stable norm as seen in section \ref{sectionnorm} and the measure part of
theorem \ref{lesfonctionsalg}.

\proof
Let $f$ be a continuous function from $\R$ to $\R$ and define 

\hfil \fonction%
f_\rho:B_\rho(1)|\R |%
x|f\bigl(d_\rho\bigl(0,x)\bigr)| \hfill

\noindent and $f_\infty(x) =f\bigl(\n{x}_\infty\bigr)$ on $B_\infty(1)$.
We want to show that (remember that $\delta_\rho(x)=\rho x$)
$$
\int f_\rho \cdot \chi_{B_g(\rho)}\circ \delta_\rho d\mu_\rho \tend_{\rho \to \infty} \int_{B_\infty(1)} f_\infty \,d\mu_\infty \numeq \label{lebut}
$$
To obtain this we are going to cut the difference in three pieces i.e. :
$$\eqalignno{
\qquad \Biggl\vert \int f_\rho \cdot \chi_{B_g(\rho)}\circ \delta_\rho~d\mu_\rho &-
\int_{B_\infty(1)} f_\infty \,d\mu_\infty \Biggr\vert \leq \hfill & \cr
  &\Biggl\vert \int f_\rho \cdot 
   \bigl(\chi_{B_g(\rho)}\circ \delta_\rho - \chi_{B_\infty(1)} \bigr)~d\mu_\rho\Biggr\vert &\nume \label{terme1}\cr
  &+\Biggl\vert 
   \int_{B_\infty(1)} \bigl(f_\rho-f_\infty \bigr)d\mu_\rho
   \Biggr\vert \qquad  \qquad \qquad \ \,& \nume \label{terme2}\cr
  &+\Biggl\vert \int_{B_\infty(1)} f_\infty~d\mu_\rho-\int_{B_\infty(1)} f_\infty~d\mu_\infty 
   \Biggr\vert & \nume \label{terme3}}
$$
Now it suffices to notice that
\item{1)} the part \ref{terme1} goes to $0$ because inside we have the product
of $\chi_{B_g(\rho)}\circ \delta_\rho - \chi_{B_\infty(1)} $, which is easily seen to simply converge to $0$ 
thanks to \ref{cor:normes}, 
 with bounded terms compactly supported.
\item{2)} Same reason but here it is $f_\rho-f_\infty$ which simply converges to $0$,
\item{3)} Finally the convergence to $0$ of \ref{terme3} is due once again to the measure
part of theorem \ref{lesfonctionsalg}.

As a conclusion we have \ref{lebut}. Injecting now $f_\rho$ into
the Raleigh's quotient we get :
$$
\rho^2 \lambda_g\bigl(B_g(\rho)\bigr) \leq\frac{\displaystyle{\int \bigl((f')_\rho\bigr)^2  \cdot \chi_{B_g(\rho)}\circ\delta_\rho~d\mu_\rho}}%
{\displaystyle{\int (f_\rho)^2  \cdot \chi_{B_g(\rho)}\circ \delta_\rho~d\mu_\rho}}
$$
we apply  \ref{lebut} to obtain
$$
\limsup_{\rho \to \infty} \rho^2 \lambda_g\bigl(B_g(\rho)\bigr) \leq \frac{\displaystyle{\int_{B_\infty(1)} \bigl((f')_\infty\bigr)^2~d\mu_\infty}}%
{\displaystyle{\int_{B_\infty(1)} f_\infty^2~d\mu_\infty}}
$$
and now taking for  $f$ the right function (i.e. the solution of the differential equation
$f''+\frac{n-1}{x}f'(x)+\lambda_{e,n}f=0$) we can conclude.

Let us now study the equality case. Take again the function $f$ which
gives the eigenfunction of the Euclidean Laplacian on the Euclidean unit ball
(i.e. the solution of $f''+\frac{n-1}{x}f'(x)+\lambda_{e,n}f=0$) and normalise it.
The $\Gamma$-convergence theory allows to say, 
 taking ${\cal E}_\rho$ and ${\cal E}_\infty$ as the energies of $\Delta_\rho$ and $\Delta_\infty$ 
on the balls $B_\rho(1)$ and  $B_\infty(1)$ respectively for the adapted measures
and thanks to proposition \ref{thdefstrspe} and \ref{theo1}
$$
{\cal E}_\infty(f_\infty)\leq \liminf_{\rho\to \infty}{\cal E}_\rho(f_\rho)\leq \limsup_{\rho\to \infty}{\cal E}_\rho(f_\rho)\leq \lambda_{e,n}\numeq \label{ana}
$$
Now from the equality assumption we have
$$
\lambda_{e,n}\leq {\cal E}_\infty(f_\infty){\rm,\ }\numeq\label{pafsh}
$$
thus \ref{ana} and \ref{pafsh} imply equality which in turn imply that 
$f_\infty$ is an eigenfunction for the first eigenvalue . 
Hence $f_\infty$ is smooth (at least in a neighbourhood of zero).

Now from the study of Bessel's function (see. {bow}, {\S}103 {\`a} {\S}105)) we see that
taking $p=(n-2)/2$ we have $f(x)=x^{-p}J_p\bigl(\sqrt{\lambda_e}x\bigr)$ with $J_p$ an analytic
function defined by (see F.~Bowman \cite{bow} {\S}84)
$$
J_p(x)=\frac{x^p}{2^p\Gamma(p+1)} \biggl( 1-\frac{x^2}{2\cdot2n+2}+\frac{x^4}{2\cdot4\cdot2n+2\cdot2n+4}+\dots\biggr)
$$ 
thus $f$ has the following shape 
$$
f(x)=\frac{\lambda_e^p}{2^p\Gamma(p+1)}\biggl( 1-\frac{x^2\lambda_e}{2\cdot2n+2}+\frac{x^4\lambda_e^2}{2\cdot4\cdot2n+2\cdot2n+4}+\dots\biggr)
$$
in other words $f$ has the following asymptotic expansion $1+\alpha_1x^2+\alpha_2x^4+\dots$ 
(up to a multiplicative constant), now notice that the function
$1+\alpha_1x+\alpha_2x^2+\dots$ admits an inverse $g \in C^\infty$  in a neighbourhood of zero,
which implies that $g\circ f_\infty(x) =cst\cdot\n{x}_\infty^2$ is $C^2$ in a neighbourhood of zero,
thus the stable norm comes from a scalar products, which means that it is
Euclidean. 

In fact we have some more informations. Indeed in order for $f_\infty$ to be an eigenfunction,
the norme of the differential of the stable norm with respect to the
albanese metric (the scalar product giving the laplacian $\Delta_\infty$) must be almost everywhere
equal to one (a simple computation using the fact that the stable norm is euclidean and
the Cauchy-Schwartz inequality). Which implies that the unit ball of the albanese
metric must be inside the unit ball of the stable norm. Now the maximum principle
and the monotony with respect to inclusion of the eigenvalues implies that
equality holds if and only if the stable norm and the Albanese metric coincides.
The stable norm being the Albanese metric we can now use the theorem \ref{sneuc} 
to conclude.
\qed

\proclaim Theorem \theo.
Let $(\T^n,g)$ be a torus, its stable norm coincides with the Albanese metric if and only if
the torus is flat.\label{sneuc}

\proof
Let us take a base $\eta_1,\dots,\eta_n$ of Harmonic one forms, any function $\alpha$
and any $2$-form $\beta$. We shall write $(\cdot,\cdot)_g$ pointwise scalar product
induced by $g$ on forms ($\n{\cdot}_g$ the associated norm) and $\langle\cdot,\cdot\rangle_g$ the integral scalar 
product normalized by the volume. Then by Hodge's theorem
$$
\n{\eta_i}_\infty^2= \inf_{\alpha,\beta} \sup_{x \in \T^n} \n{\eta_i}_g^2 + \n{d\alpha}_g^2 +\n{\delta\beta}_g^2 \geq 
\sup_{x \in \T^n} \n{\eta_i}_g^2
$$
and 
$$ 
\langle\eta_i,\eta_i\rangle_g = \frac{1}{\Vol_g(\T^n)}\int_{\T^n} \n{\eta_i}_g^2 d\vol_g \leq  \sup_{x \in \T^n} \n{\eta_i}_g^2
$$
the case of equality implies that $\langle\eta_i,\eta_i\rangle_g=(\eta_i,\eta_i)_g(x)$ for all $x\in \T^n$.
Now it suffices to see that the metric $g$ can be written in the following
way :
$$\sum_{i,j} \lambda_{ij} \eta_i\circ\eta_j = g$$
where $\eta_i\circ \eta_j = 1/2(\eta_i\otimes\eta_j+\eta_j\otimes\eta_i)$ and  
$\Lambda=(\lambda_{ij})$ is the matrice such that $\Lambda^{-1}=\bigl(\langle\eta_i,\eta_j\rangle_g\bigr)$.
Now taking local $f_i$ such that $df_i=\eta_i$ then the function $F(x)=(f_1(x),\dots,f_n(x))$
is an isometry between a open set of $\T^n$ and an eucliean space, 
thus the torus is flat.
\qed
\partie{Related topics}\label{relatedtopics}
In that section we come back to the asymptotic volume, proving in the meantime a
generalised Faber-Krahn inequality. Then we explain what can be deduces from
our work for the heat kernel and how it is related to other's work. We finally
state how theorem \ref{theo1} passes to graded nilmanifolds.

\subpartie{Asymptotic volume of tori}\label{monvolas}
\subsubpartie{Generalised Faber-Krahn inequality}\label{fbk}
We need some more definitions

\proclaim Definition \theo.
For a rectifiable submanifold $N$ of $\R^n$ (we can think of it of finite adapted Hausdorff measure)
we will write $I(N)$ the associated integral current. For an integral current $C$, $M(C)$ will be
its mass as defined par H.~Federer (see \cite{federer} for example).

\proclaim Definition \theo.
Let $\R^n$, with the norm $\n{\cdot}$ ($\n{\cdot}_*$ will be the dual norm), we define 
$$
\lambda_1\bigl(\Omega, \n{\cdot}\bigr)= \inf_{f} %
\frac{\displaystyle{\int_\Omega \n{df}^2_* d\mu}}{\displaystyle{\int_\Omega f^2 d\mu}}
$$
where $\mu$ is the Lebesgue measure on $\R^n$, and the infimum is taken
over all lipschitz functions vanishing on the border.

The following lemma holds,

\proclaim Lemma {\theo.} (Faber-Krahn inequality for norms). 
Let $D$ be a domain of $\R^n$, with the norm $\n{\cdot}$ and a measure $\mu$ invariant
by translation. Let $D^*$ be the norm's ball with same measure as $D$, then \label{lemfk}
$$
\lambda_1\bigl(D^*, \n{\cdot}\bigr) \leq %
\lambda_1\bigl(D, \n{\cdot}\bigr) \numeq \label{fk}
$$
the equality case implying that $D$ is a norm's ball.

\proof
We need two ingredients for this proof. The first is an isoperimetric inequality,
which is given by a result of Brunn (see a proof by M.~Gromov in \cite{MS}). The second
is a co-area formula, which can be found in Federer \cite{federer} p. 438.

More specifically, let us write $G_t=\{x \mid \vert f(x) \vert = t\}$ then on one side we have
$$
\int_\Omega h \alpha\land df = \int_0^{\sup f} \int_{G_t} h \alpha_{\mid G_t} dt = \int_0^{\sup f} I_{\vert f\vert =t}(h \alpha) dt
$$
and on the other
$$
\int_\Omega \n{df}_* d\mu = \int_0^{\sup f} M(I_{\vert f \vert =t}) dt \numeq \label{masse}
$$
(see P.~Pansu \cite{pansu2}) where $d\mu$ is the translation invariant volume form on $\R^n$
such that the norm's ball of radius one has measure $1$.

Take $\alpha = \displaystyle{\frac{1}{\vert df\vert^2}*df}$
where $*$ is the Hodge operator on differential forms over $\R^n$ then we get
$$
\int_\Omega h d\mu = \int_0^{\sup f} \int_{G_t} h \alpha_{\mid G_t} dt
$$
Take now a look at the same equality on $\Omega_t=\bigl\{ x \mid |f(x)|>t \bigr\}$ i.e. :
$$
\int_{\Omega_t} h d\mu = \int_t^{\sup f} \int_{G_t} h \alpha_{\mid G_t} dt%
= \int_t^{\sup f} I_{\vert f\vert =t}(h \alpha) dt\numeq \label{omegat}
$$
differentiating each member of equality \ref{omegat}
we get almost everywhere the following equality
$$
\int_{G_t} h\alpha_{\mid G_t} = I_{\vert f\vert =t}(h\alpha)\numeq\label{courant}
$$
taking into account \ref{courant} and \ref{masse} we obtain
$$
\int_{G_t} \n{df}_* \alpha_{\mid G_t} = M(I_{\vert f \vert =t})\numeq \label{blahh}
$$
applying the Cauchy-Schwartz inequality to the left side of \ref{blahh}
and make the appropriate identification thanks to \ref{courant} we finally have
$$
\frac{M(I_{\vert f \vert =t})^2}{I_{\vert f\vert=t}(\alpha)} \leq 
{I_{\vert f\vert=t}\Bigl(\bigl(\n{df}_*\bigr)^2\alpha\Bigr)}\numeq \label{fb1}
$$

The function $f^*$ associated to $f$ by symmetrisation is lipschitz,
thus she satisfies a similar co-area formula. Hence we have for almost all $t$
$$
I_{\vert f\vert=t}(\alpha) =-\frac{d}{dt} {\rm Vol} (\Omega_t) =%
 -\frac{d}{dt} {\rm Vol}(\Omega^*_t) = I_{\vert f^*\vert=t}(\alpha^*)\numeq \label{dvolu}
$$
now using the Brunn's isoperimetric inequality (see \cite{MS}) we have
$$
M(I_{\vert f^* \vert =t})\leq M(I_{\vert f \vert =t})\numeq \label{ibrun}
$$
injecting \ref{dvolu} and \ref{ibrun} in 
\ref{fb1} and noticing that $\n{df^*}_*$ is constant on $\{ \vert f \vert =t\}$,
which implies that the equivalent of \ref{fb1} for $f^*$ is an equality we get
(for almost all $t$)
$$
{I_{\vert f^*\vert=t}\Bigl(\bigl(\n{df^*}_*\bigr)^2\alpha^*\Bigr)} =
\frac{M(I_{\vert f^* \vert =t})^2}{I_{\vert f^*\vert=t}(\alpha^*)} \leq
\frac{M(I_{\vert f \vert =t})^2}{I_{\vert f\vert=t}(\alpha)} \leq 
{I_{\vert f\vert=t}\Bigl(\bigl(\n{df}_*\bigr)^2\alpha\Bigr)}\numeq \label{fbouf}
$$
Now we sum the extremal terms of \ref{fbouf} to obtain the wanted
inequality :
$$
\int_{\Omega^*} \bigl(\n{df^*}_*\bigr)^2 dv \leq %
\int_\Omega \bigl(\n{df}_*\bigr)^2 dv
$$
which allows us to conclude the proof because
$$\int_{\Omega^*} (f^*)^2 dv = \int_\Omega (f)^2 dv$$
For the equality case, it suffices to see that it implies the equality
case in the Brunn's isoperimetric inequality to conclude.
\qed
Let us notice that this lemma immediately implies
\proclaim Corollary.
  Let $D_1$ the unit ball of the norm $\n{\cdot}$ then
$$ \lambda_1(D_1)=\lambda_{e,n}$$ thus
$$
\lambda_{e,n} \biggl(\frac{\mu(D_1)}{\mu(D)}\biggr)^{\frac{2}{n}} \leq %
\lambda_1\bigl(D, \n{\cdot}\bigr)
$$
where $\mu$ is a Haar measure on $\R^n$.\label{corlemfk}

\proof
The symmetrisation from the previous theorem shows that the minimum of the 
Rayleigh's quotient
is obtained with functions depending on the distance from the centre of the ball. Hence we
are lead to the same calculations as in the Euclidean case.
\qed

\subsubpartie{Lower bound for the asymptotic volume}\label{asvol}

We are now going to apply the generalised Faber-Krahn inequality to $\lambda_\infty$.
In that aim in mind let us notice that $\lambda_\infty\bigl(B_\infty(1)\bigr)=\lambda_1\bigl(B_\infty(1),\n{\cdot}_2^*\bigr)$
with the dual norm of $\n{\cdot}_2^*$ defined by
$$
\n{\xi}_2=\sum_{ij} q^{ij}\xi_i\xi_j
$$
and let us write $B_{\rm Al}$ the unit ball of $\n{\cdot}_2^*$.
We now can apply inequality of lemma \ref{lemfk} and more precisely its corollary \ref{corlemfk}
$$
\lambda_\infty\bigl(B_\infty(1)\bigr)\geq \left(\frac{\mu(B_{\rm Al})}{\mu\bigl(B_\infty(1)\bigr)}\right)^{2/n} \lambda_{e,n}
$$
where $\mu$ is any Haar measure. Now applying theorem \ref{lambdaas} we get
$$
\left(\frac{\mu(B_{\rm Al})}{\mu\bigl(B_\infty(1)\bigr)}\right)^{2/n} \lambda_{e,n} \leq  \lambda_{e,n} \numeq\label{backtovol}
$$

We finally get the following proposition taking in \ref{backtovol}  the Haar measure such that
the measure of $B_\infty(1)$ is the asymptotic volume (i.e. the measure $\mu_\infty$) and transforming the other term
in order to make the Albanese's Torus's volume appear.

\proclaim Proposition {\ref{burivame}}.
Let $(\T^n,g)$ be a Riemannian Torus, $B_g(\rho)$ the geodesic balls of radius $\rho$
centred on a fixed point and ${\rm Vol}_g\bigl(B_g(\rho)\bigr)$ their Riemannian volume 
induced on the universal cover, writing
$$
{\rm Asvol}(g) =\lim_{\rho\to \infty} \frac{{\rm Vol}_g\bigl(B_g(\rho)\bigr)}{\rho^n}
$$
then
\item{1.} $ {\rm Asvol}(g) \geq \frac{\displaystyle{{\rm Vol}_g(\T^n)}}{\displaystyle{{\rm Vol}_{{\rm Al}}(\T^n)}}\omega_n $
\item{2.} In case of equality, the torus is flat.
Where $\omega_n$ is the unit Euclidean ball's Euclidean volume.

\proof
There remains the equality case to be proved, which can be obtained using either the
equality case of the Faber-Krahn inequality, which says that $B_\infty(1)$ is an
ellipsoid either the equality case of theorem \ref{theo2} and then
we conclude by using theorem \ref{sneuc}.
\qed

We still have two remarks concerning this proposition, the first one is
included in the following corollary.

\proclaim Corollary.
  For $n=2$ we have
\item{1.} $ {\rm Asvol}(g) \geq \pi=\omega_2 $
\item{2.} In case of equality, the torus is flat.

In other words we obtain the theorem of D.~Burago et S.~Ivanov on the asymptotic volume
of tori in the $2$ dimensional case (see \cite{bi2}). The second remark is that
we can not do better this way. See \cite{vernicos} part three for more details.

\subpartie{Long time asymptotics of the heat kernel}\label{lgtime}

Let $(\T^n,g)$ be a torus and $(\R^n,\tilde{g})$ its universal cover with the lifted metric.
We remind the reader that $g_\rho=(1/\rho^2)\delta_\rho^*\tilde{g}$ are the rescaled metrics and $\Delta_\rho$
their Laplacian, here it will be on $\R^n$.

We are going to study from the homogenisation point of view the long time asymptotic
behaviour of the heat kernel i.e. we are interested in the behaviour as $t$ goes to
infinity of a solution $u(t,x)$ of the following problem
$$
\left\{\eqalign{
{\partial u\over \partial t}+ \Delta u = 0 &\qquad  {\rm in \ } \mathopen]0,+\infty\mathclose[ \times\R^n \cr
u(0,x) = u_0(x) &\cr}\right. \numeq \label{eq:chaleurA1}
$$
For a probabilistic insight one could see M. Kotani et T. Sunada \cite{kosu}.

Let us introduce the rescaled functions
$$
u_\rho(t,x)=\rho^nu(\rho^2t,\delta_\rho x), \ \rho>0
$$
it is the straightforward that (see. \ref{respell}) 
$u$ is a solution of  \ref{eq:chaleurA1} if and only if 
$u_\rho$ is a  solution of
$$
\left\{ 
\eqalign{%
{\partial u_\rho\over \partial t} + \Delta_\rho u_\rho = 0 & {\rm in\ } \mathopen]0,+\infty\mathclose[ \times\R^n \cr
u_\rho(0,x) = \rho^nu_0(\delta_\rho x) & \cr }\right. \numeq \label{eq:chaleurAr}
$$
hence studying $u(t,\cdot)$ as $t$ goes to infinity is the same as
studying $u_\rho(1,\cdot)$ as $\rho\to \infty$. In other words we are once again lead to the
study of the spectral structures $(\Delta_\rho)$ on $\R^n$. We have
\proclaim Theorem \theo.
  The net of resolvents $(R_\lambda^\rho)$ weakly converges to the
resolvent $(R_\lambda^\infty)$ of $\Delta_\infty$ in $L^2(\R^n)$.

\noindent{\bf Remark.} The proof is the same as \ref{cvres}. In fact in
that case we would rather talk of {\sl G-convergence}.
We now can apply the theorems from the chapter III of \cite{gconv}
more precisely theorems $4$ and $6$.

\proclaim Theorem {\theo.}(\cite{gconv} page 136).
 The fundamental solution  $k(t,x,y)$  of \ref{eq:chaleurA1} has he following
asymptotic expansion
$$
k(t,x,y)=k_\infty(t,x,y)+t^{-\frac{n}{2}} \theta(t,x,y)
$$ 
where $k_\infty(t,x,y)$ is fundamental solution of
$$
{\partial u_\infty\over \partial t} + \Delta_\infty u_\infty = 0  {\rm\ in \ } \mathopen]0,+\infty\mathclose[ \times\R^n \numeq \label{eq:chaleurAinf}
$$
and $\theta(t,x,y)\to 0$ uniformly as $t\to \infty$ on $|x|^2+|y|^2\leq at$, for any
fixed constant $a>0$.
      
\noindent{\bf Remark.}
This is slightly weaker that theorem~1 of M.~Kotani
and T.~Sunada in \cite{kosu}.

\proclaim Theorem {\theo.} (\cite{gconv} page 138).
Let $u_0 \in L^1(\R^n)\cap L^\infty(\R^n)$. Then  $u(t,x)$ the solution of\ref{eq:chaleurA1} 
has the following asymptotic expansion :
$$
  u(t,x) = c_0 (4\pi t)^{-\frac{n}{2}} \int_{\R^n} u_0(y) dy + t^{-\frac{n}{2}} \theta(t,x)
$$
where $\theta(t,x)$ converges uniformly to $0$ for $|x|<R$ where $R$ is
a positive constant and $c_0$ is the determinant of the matrix associated to $\Delta_\infty$.

That last claim can be precised by
\proclaim Theorem {\theo.} (Duro, Zuazua \cite{duzu}).
Let $u_0\in L^1(\R^n)$. The sole solution of
\ref{eq:chaleurA1} satisfies for every $p\in [1,+\infty[$ :
$$
  t^{n/2(1-1/p)}\n{u(t)-u_\infty(t)}_p \to 0,\ {\rm as\  } t \to +\infty \numeq \label{eq:chaleur:homo1}
$$
where $u_\infty$ is the unique solution of the homogenised problem \ref{eq:chaleurAinf}. For
$n=1$ and $n=2$ \ref{eq:chaleur:homo1} is also true for $p=\infty$.

\subpartie{The macroscopical sound of graded nilmanifolds}\label{grnilm}
In this part we want to emphasise the fact that theorem \ref{theo1} is still
true for graded nilmanifolds, at least for the Dirichlet case, but it involves some
sub-riemannian geometry. We just give the statement the details are
to be found in \cite{vernicos} chapter two.

\proclaim Theorem \theo.
Let $(M^n,g)$ be graded nilmanifold, $B_g(\rho)$ the induced Riemannian
ball of radius $\rho$
on its universal cover and $\lambda_i\bigl(B_g(\rho)\bigr)$ the $i^{\rm th}$
eigenvalue of the Laplacian on $B_g(\rho)$ for the Dirichlet problem.
\endgraf
Then there exists an hypoelliptic operator $\Delta_\infty$ (the Kohn Laplacian
of a left invariant metric), whose
$i^{\rm th}$ eigenvalue for the Dirichlet problem on the stable ball
is $\lambda_i^\infty$ and such that
$$
\lim_{\rho\to \infty} \rho^2 \lambda_i\bigl(B_g(\rho)\bigr)= \lambda_i^\infty
$$

\noindent where the stable ball is the metric ball given by 
the Carnot-Caratheodory distance
found in \cite{pansu} and arising from the stable norm.

\partie{Proof of theorem \ref{cvssets}}\label{lastbutnotleast}

Let $(\Sigma_\alpha)$ be a net of spectral structures and let us focus on the spectra. For
a fixed operator $\sigma(\cdot)$ will be its spectrum. Let us begin with the case
of strong convergence.
\proclaim Proposition \theo.
  If $\Sigma_\alpha \to \Sigma$ strongly, then for any $\lambda\in \sigma(A)$ there is $\lambda_\alpha\in \sigma(A_\alpha)$ such that
the net $(\lambda_\alpha)$ converges to $\lambda$, this is written 
$$
\sigma(A)\subset\lim_\alpha \sigma(A_\alpha)
$$

\proof
  Let $\lambda\in \sigma(A)$ and $\varepsilon>0$ and take $\zeta=\lambda+i\varepsilon$ then : 
$$\n{R^\alpha_\zeta}_{{\cal L}_\alpha} = \frac{1}{\inf_{\rho\in \sigma(A_\alpha)} |\zeta-\rho|} \quad {\rm and } \quad%
\n{R_\zeta}_{{\cal L}_\infty}=\frac{1}{\inf_{\rho\in \sigma(A)} |\zeta-\rho|}=\frac{1}{\varepsilon}{\rm .}$$
From the assumption, the net of resolvents strongly converges hence by \ref{ncvsnf}
$$
\limsup_\alpha \inf_{\rho \in \sigma(A_\alpha)} |\zeta-\rho| \leq \varepsilon
$$
and as it is true for any $\varepsilon$, we can conclude.
\qed

\proclaim Lemma \theo.
For any reals $a$,$b$ out of the spectra of $A$ such that $-\infty\leq a<b\leq+\infty$ then \label{leseize}
$$
a\leq\frac{{\cal E}(u)}{\n{u}^2_\infty}\leq b \quad {\sl pour\ tout\ } u\in E\bigl(\mathopen]a,b]\bigr)L_\infty^2\setminus\{0\}.
$$
(where $E\bigl(\mathopen]a,b]\bigr)=E\bigl(\mathopen]a,+\infty\mathclose[\bigr)$ si $b=+\infty$).

\proof
Let $a<b$ two reals out of the spectra of $A$ and
  $$ u\in E\bigl(\mathopen]a,b]\bigr)L_\infty^2\setminus\{0\}{\rm .}$$
then
$$
\int_{\mathopen]a,b]}dEu= E\bigl(\mathopen]a,b]\bigr)u=u=\int_\R dEu
$$
thus $\scal{Eu,u}=0$ on $\R\setminus\mathopen]a,b]$. Now if $u\in D(A)$,
$$
{\cal E}(u)=\scal{Au,u}=\int_\R \lambda\ d\scal{E(\lambda)u,u}=\int_{\mathopen]a,b]} \lambda\ d\scal{E(\lambda)u,u}
$$
and the last term satisfies
$$
a\n{u}^2_\infty=a\int_{\mathopen]a,b]} d\scal{E(\lambda)u,u}\leq\int_{\mathopen]a,b]} \lambda\ d\scal{E(\lambda)u,u}\leq b\int_{\mathopen]a,b]} d\scal{E(\lambda)u,u}=b\n{u}^2_\infty
$$
\qed

For any Borel set $I\subset \R$ we write $n(I)=\dim E(I)L_\infty^2$ and $n_\alpha(I)=\dim E_\alpha(I)L^2_\alpha$.
\proclaim Proposition \theo.
  Let $a<b$ two reals out of the point spectrum of $A$. if $\Sigma_\alpha\to \Sigma$ strongly then
$$
\liminf_\alpha n_\alpha\bigl(\mathopen]a,b]\bigr) \geq n\bigl(\mathopen]a,b]\bigr)
$$
 and in particular,\label{ladeuxsix}
$$
\liminf_\alpha \dim L^2_\alpha\geq \dim L^2_\infty
$$

\proof
Lets take an orthonormal base $\{ \varphi_k \mid k=1,\dots,n\bigl(\mathopen]a,b]\bigr)\}$  
of $E\bigl(\mathopen]a,b]\bigr)L^2_\infty$.
Let $n\in \N$ a fixed number if  $n\bigl(\mathopen]a,b]\bigr)=\infty$ else $n=n\bigl(\mathopen]a,b]\bigr)$.
Then there are nets $\varphi_k^\alpha\in L^2_\alpha$ for $k=1,\dots,n$ such that $\lim_\alpha \varphi_k^\alpha=\varphi_k$.
As $E_\alpha\bigl(\mathopen]a,b]\bigr) \to E\bigl(\mathopen]a,b]\bigr)$ strongly, taking
$\psi_k^\alpha=E_\alpha\bigl(\mathopen]a,b]\bigr)\varphi_k^\alpha$ we get
$$
\lim_\alpha \psi_k^\alpha=E\bigl(\mathopen]a,b]\bigr)\varphi_k=\varphi_k
$$
hence
$$
\lim_\alpha\scal{\psi_i^\alpha,\psi_j^\alpha}_\alpha=\scal{\varphi_i,\varphi_j}=\delta_{ij}
$$
from which we deduce that $(\psi_k^\alpha)_{k=1,\dots,n}$ is a free family for $\alpha$ large enough and
$$
\lim_\alpha n_\alpha\bigl(\mathopen]a,b]\bigr)\geq n.
$$
which proves the first assertion.
For the second it comes from the fact that
$n\bigl(\mathopen]a,b]\bigr)$ converges to $\dim L^2_\infty$ as $a\to -\infty$ and $b\to +\infty$.
\qed
Let us now have a look at the compact convergence case :

\proclaim Theorem \theo.
  If $\Sigma_\alpha \to \Sigma$ compactly converges, then for any $a$,$b$ out of the point spectrum
of $A$ such that $a<b$ for  $\alpha$ large enough we have
$n_\alpha\bigl(\mathopen]a,b]\bigr)=n\bigl(\mathopen]a,b]\bigr)$. In  particular the limit of the sets 
$\sigma(A_\alpha)$ coincides with $\sigma(A)$ \label{lespec}

\proof
The compact convergence implies that the operators $R_\zeta$, $T_t$ and $E\bigl(\mathopen]\lambda,\mu]\bigr)$ are
compact (see \ref{ncvsnf}) 
Thus the spectrum of $A$ is discrete and $n\bigl(\mathopen]a,b]\bigr)<\infty$ if $a<b<\infty$.
Let $(0\leq)\lambda_1\leq\lambda_2\leq\dots \leq\lambda_n$ be  the spectrum of $A$, where 
$$
\cases{
  n=  0 & if the spectre is empty , \cr
  n \in \N & if the spectre is finite and\cr
  n=\infty & if the spectre is a sequence converging to infinity.}
$$

{\bf Step 1 : }
Fix $\varepsilon_0$ and let $\Lambda_1^\alpha=E\bigl(\mathopen]-\infty,\lambda_1+\varepsilon_0]\bigr)L^2_\alpha$ and $\Lambda_1=L^2_\infty$,
where \hbox{$\lambda_1=\lambda_1+\varepsilon_0=\infty$} if $n=0$. Let
$$
\mu_1= \liminf_\alpha \inf \bigl\{ {\cal E}_a(u) \mid \n{u}_\alpha=1,\ u\in \Lambda_1^\alpha \bigr\}
$$
lemma \ref{leseize} allows us to say that $\lim_\alpha n_\alpha\bigl(\mathopen]-\infty,\mu]\bigr)=0$ for any 
$\mu\in\mathopen]-\infty,\mu_1\mathclose[$. Applying proposition \ref{ladeuxsix} we get
$n\bigl(\mathopen]-\infty,\mu]\bigr)=0$, in other words
for any $\mu\leq \mu_1$ then $\mu\leq \lambda_1$ thus $\mu_1\leq \lambda_1$. 
Hence if $\mu_1=+\infty$, $n=0$ and $L^2_\alpha={0}$ for $\alpha$ large enough and the theorem is proved
in that case.

Suppose that $\mu_1<+\infty$. For $\alpha$ large enough
we can find unit vectors $\varphi_1^\alpha\in \Lambda_1^\alpha$ such that $\liminf_\alpha {\cal E}_\alpha(\varphi_1^\alpha)=\mu_1$. From the asymptotic
compactness of ${\cal E}_\alpha$ we can extract a sub-net $(\varphi_1^\alpha)_{\alpha\in{\cal A}}$  such that $\varphi_1=\lim_\alpha \varphi_1^\alpha$
strongly and thanks to definition \ref{thdefstrspe} ${\cal E}(\varphi_1)\leq\mu_1$.
The strong convergence induces the convergence of the norms hence $\n{\varphi_1}= 1$ and 
$$
\lambda_1=  \inf \bigl\{ {\cal E}(u) \mid \n{u}=1,\ u\in \Lambda_1 \bigr\} \leq {\cal E}(\varphi_1)\leq\mu_1< +\infty.
$$
As a consequence $n\geq1$, $\lambda_1=\mu_1={\cal E}(\varphi_1)$ and $\varphi_1$ is eigenvector of $A$ for $\lambda_1$.

Furthermore let us notice that as $E_\alpha\bigl( \mathopen]\lambda_1-\epsilon,\lambda_1+\epsilon] \bigr) \to E\bigl( \mathopen]\lambda_1-\epsilon,\lambda_1+\epsilon] \bigr)$
strongly for any $\epsilon>0$ fixed and  $E\bigl( \mathopen]\lambda_1-\epsilon,\lambda_1+\epsilon] \bigr) \to E\bigl(\{\lambda_1\}\bigr)$ strongly
when $\epsilon\to 0$
there is a net of positives numbers $\epsilon_1^\alpha\to 0$ such that
$E_\alpha\bigl( \mathopen]\lambda_1-\epsilon_1^\alpha,\lambda_1+\epsilon_1^\alpha] \bigr) \to E\bigl(\{\lambda_1\}\bigr)$ strongly. From this we obtain
a net
$$
  \psi_1^\alpha=E_\alpha\bigl( \mathopen]\lambda_1-\epsilon_1^\alpha,\lambda_1+\epsilon_1^\alpha] \bigr)\varphi_1^\alpha \to E\bigl(\{\lambda_1\}\bigr)\varphi_1=\varphi_1.
$$

{\bf Step 2 :}
Let $\Lambda_2^\alpha=E\bigl(\mathopen]-\infty,\lambda_2+\varepsilon_0]\bigr)L^2_\alpha\cap \langle\varphi_1^\alpha\rangle^\bot$, $\Lambda_2= \langle\varphi_1\rangle^\bot$ and
$$
\mu_2= \liminf_\alpha \inf \bigl\{ {\cal E}_a(u) \mid \n{u}_\alpha=1,\ u\in \Lambda_2^\alpha \bigr\}.
$$
Again lemma \ref{leseize} allows us to say that $\lim_\alpha n_\alpha\bigl(\mathopen]-\infty,\mu]\bigr)=0$
for any $\mu \in \bigl(\mathopen]\mu_1,\mu_2 \mathclose[\bigr)$  and proposition \ref{ladeuxsix} that $\mu_2\leq \lambda_2$.
Hence if $\mu_2=+\infty$, we have $n=1$ and $L^2_\alpha=\langle\psi_1^\alpha\rangle$ for $\alpha$ large enough. 
Assume $\mu_2<\infty$. Take the unitary vectors $\varphi_2^\alpha\in \Lambda_2^\alpha$  such that $\liminf_\alpha {\cal E}_\alpha(\varphi_2^\alpha)=\mu_2$. 
Then the same discussion as step 1 gives $n\geq2$, $\lambda_2=\mu_2$  and the strong convergence of a sub-net
of $(\varphi_2^\alpha)$ to $\varphi_2$ an eigenvector of $A$ for the eigenvalue $\lambda_2$. We also find a net $\epsilon_2^\alpha\to 0$
such that $\psi_2^\alpha=E_\alpha\bigl( \mathopen]\lambda_2-\epsilon_2^\alpha,\lambda_2+\epsilon_2^\alpha] \bigr)L^2_\alpha \to  \varphi_2$. Now let us notice that
for any $\epsilon>0$ there is $\alpha_\epsilon\in {\cal A}$ such that for all $\alpha\geq\alpha_\epsilon$ we have
    \item{(1)} $\psi_i^\alpha\in E_\alpha\bigl( \mathopen]\lambda_i-\epsilon,\lambda_i+\epsilon] \bigr)L^2_\alpha$ for $i=1$,$2$ ;
    \item{(2)} if $\lambda_1+2\epsilon< \lambda_2$ then
$$
E_\alpha\bigl( \mathopen]\lambda_1-\epsilon,\lambda_1+\epsilon]\bigr)L^2_\alpha = \langle\psi_1^\alpha\rangle \quad {\rm and} %
\quad E_\alpha\bigl( \mathopen]\lambda_1+\epsilon,\lambda_2-\epsilon]\bigr)L^2_\alpha ={0}. 
$$

{\bf Step 3 :} We repeat this procedure. Setting
$$
\Lambda_k^\alpha=E\bigl(\mathopen]-\infty,\lambda_k+\varepsilon_0]\bigr)L^2_\alpha\cap \langle\psi_1^\alpha,\dots,\psi_{k-1}^\alpha\rangle^\bot
$$
we have
$$
\lambda_k=\mu_k=\liminf_\alpha \inf \bigl\{ {\cal E}_a(u) \mid \n{u}_\alpha=1,\ u\in \Lambda_k^\alpha \bigr\}
$$
for $k\leq n$. Let $k\in \{ 1,2,\dots,n\}$ and $\epsilon>0$ be sufficiently small compared with $k$.
Then, there exists $\alpha_{k,\epsilon}\in {\cal A}$ such that for any $\alpha\geq \alpha_{k,\epsilon}$,
    \item{(1)} for each $\lambda \in \{\lambda_1,\dots,\lambda_{k-1}\}$ with $\lambda< \lambda_k$,
  $$
  E_\alpha\bigl( \mathopen]\lambda-\epsilon,\lambda+\epsilon] \bigr)L^2_\alpha= \langle\psi_i^\alpha \mid p_\lambda\leq i\leq q_\lambda\rangle,
  $$
  where $p_\lambda=\min\{ i\in \N \mid \lambda_i=\lambda\} $ and $q_\lambda=\max\{ i\in \N \mid \lambda_i=\lambda\}$ ;
    \item{(2)} for each $i=1,\dots,k-1$ with $\lambda_i<\lambda_{i+1}$,
$$
E_\alpha\bigl( \mathopen]\lambda_i+\epsilon,\lambda_{i+1}-\epsilon] \bigr)L^2_\alpha=\{0\}. 
$$

{\bf Conclusion} Let $a$,$b \in \R^+\setminus\sigma(A)$ two given real numbers such that $a<b$, then from what
precedes we have for $\alpha$ large enough
$$
E_\alpha\bigl( \mathopen]a,b] \bigr)L^2_\alpha= \langle\psi_k^\alpha \mid k=1,\ldots,n {\rm\ with\ } a<\lambda_k\leq b\rangle.
$$
Thus $n_\alpha\bigl( \mathopen]a,b] \bigr)$ coincides with the number $k$ such that $a<\lambda_k\leq b$,
in other words $n\bigl( \mathopen]a,b] \bigr)$.
\qed

The proof of theorem \ref{cvssets} is the same as above, but defining the $\Lambda_k^\alpha$ with the
help of the $\varphi_k^\alpha$. 
\bibliography


\bibitem[BD98]{BrDef}
\bgroup\bf A.~Braides\egroup{} and \bgroup\bf A.~Defranceschi\egroup{}.
\newblock {\em Homogenization of Multiple Integrals}.
\newblock Oxford Science Publications, 1998.
\finitem

\bibitem[Ber86]{berard}
\bgroup\bf P.~H. Berard\egroup{}.
\newblock {\em Spectral Geometry: Direct and Inverse Problems}.
\newblock Number 1207 in Lecture Notes Math. Springer Verlag, 1986.
\finitem

\bibitem[BI94]{bi1}
\bgroup\bf D.~Burago\egroup{} and \bgroup\bf S.~Ivanov\egroup{}.
\newblock Riemannian tori without conjugate points are flat.
\newblock {\em GAFA}, 4(3):259--269, 1994.
\finitem

\bibitem[BI95]{bi2}
\bgroup\bf D.~Burago\egroup{} and \bgroup\bf S.~Ivanov\egroup{}.
\newblock On asymptotic volume of tori.
\newblock {\em GAFA}, 5(5):800--808, 1995.
\finitem

\bibitem[BLP78]{blp}
\bgroup\bf A.~Bensousan\egroup{}, \bgroup\bf J.-L. Lions\egroup{}, and
  \bgroup\bf G.~Papanicolaou\egroup{}.
\newblock {\em Asymptotic analysis for periodic structures}.
\newblock Studies in mathematics and its applications. North Holland, 1978.
\finitem

\bibitem[Bow58]{bow}
\bgroup\bf F.~Bowman\egroup{}.
\newblock {\em Introduction to Bessel Functions}.
\newblock Dover, 1958.
\finitem


\bibitem[Bur92]{bu1}
\bgroup\bf D.~Burago\egroup{}.
\newblock Periodic metrics.
\newblock {\em Advances in soviet mathematics}, 9:205--210, 1992.
\finitem


\bibitem[Cha94]{chav1}
\bgroup\bf I.~Chavel\egroup{}.
\newblock {\em Riemannian geometry : a modern introduction}.
\newblock Cambridge tracts in mathematics. Cambridge University Press, 1994.
\finitem

\bibitem[DZ00]{duzu}
\bgroup\bf G.~Duro\egroup{} and \bgroup\bf E.~Zuazua\egroup{}.
\newblock Large time behavior for convection-diffusion equations in $\R^n$ with
  periodic coefficients.
\newblock {\em Journal of Differential Equation}, 167:275--315, 2000.
\finitem

\bibitem[Fed69]{federer}
\bgroup\bf H.~Federer\egroup{}.
\newblock {\em Geometric Measure Theory}.
\newblock Springer Verlag, 1969.
\finitem


\bibitem[Fuk90]{fukaya}
\bgroup\bf Fukaya\egroup{}.
\newblock Haussdorff convergence of Riemann manifold and its application.
\newblock In \bgroup\bf T.~Ochiai\egroup{}, editor, {\em Recent topic in
  Differential and Analytic Geometry}, volume 18-I of {\em Advanced Studies in
  Pure Mathematics}, pages 143--238. 1990.
\finitem

\bibitem[GHL90]{ghl}
\bgroup\bf S.~Gallot\egroup{}, \bgroup\bf D.~Hulin\egroup{}, and \bgroup\bf
  J.~Lafontaine\egroup{}.
\newblock {\em Riemannian geometry}.
\newblock Universitext. Springer-Verlag, second edition, 1990.
\finitem

\bibitem[KKO97]{kakuo}
\bgroup\bf A.~Kasue\egroup{}, \bgroup\bf H.~Kumura\egroup{}, and \bgroup\bf
  Y.~Ogura\egroup{}.
\newblock Convergence of heat kernels on a compact manifold.
\newblock {\em Kyushu J. Math}, 51:453--524, 1997.
\finitem

\bibitem[KS]{kushi}
\bgroup\bf K.~Kuwae\egroup{} and \bgroup\bf T.~Shioya\egroup{}.
\newblock Convergence of spectral structures: a functional analytic theory and
  its applications to spectral geometry.
\newblock prepublication.
\finitem

\bibitem[KS00]{kosu}
\bgroup\bf M.~Kotani\egroup{} and \bgroup\bf T.~Sunada\egroup{}.
\newblock Albanese maps and off diagonal long time asymptotics for the heat
  kernel.
\newblock {\em Commun. Math. Phys.}, 209:633--670, 2000.
\finitem

\bibitem[Laf74]{jaklaf}
\bgroup\bf J.~Lafontaine\egroup{}.
\newblock Sur le volume de la vari{\'e}t{\'e} de jacobi d'une vari{\'e}t{\'e}
  riemannienne.
\newblock {\em C. R. Acad. Sc. Paris}, 278:1519--1522, 1974.
\finitem

\bibitem[Mas93]{Dmaso}
\bgroup\bf Dal Maso\egroup{}.
\newblock {\em An Introduction to $\Gamma$-convergence.}
\newblock Birkh{\"a}user, 1993.
\finitem


\bibitem[Mos94]{mosco}
\bgroup\bf U.~Mosco\egroup{}.
\newblock Composite media and asymptotic dirichlet forms.
\newblock {\em J. Funct. Anal.}, 123(2):368--421, 1994.
\finitem

\bibitem[MS86]{MS}
\bgroup\bf V.~D. Milman\egroup{} and \bgroup\bf G.~Schechtman\egroup{}.
\newblock {\em Asymptotic Theory of finite dimensional Normed Spaces}.
\newblock Number 1200 in Lecture Notes Math. Springer Verlag, 1986.
\finitem

\bibitem[Pan82]{pansu}
\bgroup\bf P.~Pansu\egroup{}.
\newblock {\em G{\'e}ometrie du groupe de Heisenberg}.
\newblock Th{\`e}se de docteur 3{\`e}me cycle, Universit{\'e} Paris VII, 1982.
\finitem

\bibitem[Pan99]{pansu2}
\bgroup\bf P.~Pansu\egroup{}.
\newblock Profil isop{\'e}rim{\'e}trique, m{\'e}triques p{\'e}riodiques et
  formes d'{\'e}qui\-libre des cristaux.
\newblock pr{\'e}publication d'orsay, 1999.
\finitem

\bibitem[Rud91]{rudinfa}
\bgroup\bf W.~Rudin\egroup{}.
\newblock {\em Functional Analysis}.
\newblock International series in Pure and Applied Mathematics. Mc Graw-Hill,
  second edition, 1991.
\finitem


\bibitem[Ver01]{vernicos}
\bgroup\bf C.~Vernicos\egroup{}.
\newblock {\em Spectres asymptotiques des nilvari{\'e}t{\'e}s gradu{\'e}es}.
\newblock Th{\`e}se de doctorat, Universit{\'e} Grenoble I, Joseph Fourier,
  2001.
\finitem

\bibitem[ZKON79]{gconv}
\bgroup\bf V.V. Zhikov\egroup{}, \bgroup\bf S.M. Kozlov\egroup{}, \bgroup\bf
  O.A. Oleinik\egroup{}, and \bgroup\bf Kha~T'en Ngoan\egroup{}.
\newblock Averaging and g-convergence of differential operators.
\newblock {\em Russian Math. Surveys}, 34(5):69--147, 1979.
\finitem
\bibliend

\bigskip

\vbox{\hsize 7cm  \parindent=0cm \parskip=0cm
{\bf Constantin Vernicos}\par
\sf Universit\'e de Neuch\^atel\par
Institut de Math\'ematiques\par
11, rue Emile Argand\par
2007 Neuch\^atel\par
Switzerland\par
mail: \tt Constantin.Vernicos@unine.ch}

\bye

%% file: Macros.tex

\magnification 1200
\parskip 9pt plus 5pt minus 3pt
\def\item{\parskip 5pt plus 3pt minus 3pt \par\hang\textindent}

\headline{\ifodd\pageno\ifnum\pageno=1\else\hfil\rlap{\hautdepage\folio}\fi\else\llap{\hautdepage\folio}\hfil\fi} 
\footline{}              
\catcode`\@=11
\newdimen\margereliured \margereliured=1truecm
\newdimen\margereliureg \margereliureg=-0.1truecm
\def\plainoutput{%
\shipout\vbox{
\ifodd\pageno \moveright\margereliured
              \else \moveleft\margereliureg \fi
\hbox{\vbox{\makeheadline\pagebody\makefootline}}}%
\advancepageno
\ifnum\outputpenalty>-\@MM \else\dosupereject\fi}
\catcode`\@=12

\def\tenpoint{%
  \textfont0=\tenrm \scriptfont0=\sevenrm \scriptscriptfont0=\fiverm
  \def\rm{\fam0\tenrm}%
  \textfont1=\teni \scriptfont1=\seveni \scriptscriptfont1=\fivei
  \def\oldstyle{\fam1\teni}%
  \textfont2=\tensy \scriptfont2=\sevensy \scriptscriptfont2=\fivesy
  \textfont\itfam=\tenit
  \def\it{\fam\itfam\tenit}%
  \def\sl{\fam\slfam\tensl}%
  \textfont\slfam=\tensl \scriptfont\slfam=\sevensl \scriptscriptfont\slfam=\fivesl
  \def\bf{\fam\bffam\tenbf}%
  \textfont\bffam=\tenbf \scriptfont\bffam=\sevenbf
  \scriptscriptfont\bffam=\fivebf  
  \def\tt{\fam\ttfam\tentt}%
  \textfont\ttfam=\tentt
  \abovedisplayskip=6pt plus 2pt minus 6pt
  \abovedisplayshortskip=0pt plus 3pt
  \belowdisplayskip=6pt plus 2pt minus 6pt
  \belowdisplayshortskip=7pt plus 3pt minus 4pt
  \smallskipamount=3pt plus 1pt minus 1pt
  \medskipamount=6pt plus 2pt minus 2pt
  \bigskipamount=12pt plus 4pt minus 4pt
  \normalbaselineskip=12pt
  \setbox\strutbox=\hbox{\vrule height8.5pt depth3.5pt width0pt}%
  \normalbaselines\rm}
\catcode`\|=13
\def\today{\ifcase\month\or january \or february \or march \or april
\or may \or june\or july\or august \or september\or october\or november\or
December\fi\ \number\day , \number\year}


\newskip\LastSkip
\def\nobreakatskip{\relax\ifhmode\ifdim\lastskip>0pt
  \LastSkip\lastskip\unskip
  \nobreak\hskip\LastSkip
  \fi\fi}
\catcode`\;=\active \def;{\nobreakatskip\string;}
\catcode`\:=\active \def:{\nobreakatskip\string:}
\catcode`\!=\active \def!{\nobreakatskip\string!}
\catcode`\?=\active \def?{\nobreakatskip\string?}

\newif\ifrefvis
\refvisfalse

\newif\ifarxiv
\arxivfalse
\def\arXiv{\arxivtrue}

\newif\iffrance

\def\anglais{\francefalse}

\newskip\afterskip
\catcode`\@=11
\def\p@int{.\par\vskip\afterskip\penalty100} 
\def\p@intir{\discretionary{.}{}{.\kern.35em---\kern.7em}}
\def\pointir{\afterassignment\pointir@\global\let\next=}
\def\pointir@{\ifx\next\par\p@int\else\p@intir\fi\egroup\next}
\catcode`\@=12
\def|{\relax\ifmmode\vert\else\findef\fi}
\def\findef{\errhelp{Cette barre verticale ne correspond ni a un \vert mathematique
                        ni a une fin de definition, le contexte doit vous indiquer ce qui manque.
                        Si vous vouliez inserer un long tiret, le codage recommande est ---,
                        dans tous les cas, la barre fautive a ete supprimee.}%
                        \errmessage{Une barre verticale a ete trouvee en mode texte}}

\def\TITR#1|{\null{\mss\baselineskip=17pt
                           \vskip 3.25ex plus 1ex minus .2ex
                           \leftskip=0pt plus \hsize
                           \rightskip=\leftskip
                           \parfillskip=0pt
                           \noindent #1
                           \par\vskip 2.3ex plus .2ex}}
 
\def\auteur#1|{\penalty 500
               \vbox{\centerline{
                 \iffrance par \else by \fi #1}
                \vskip 10pt}\penalty 500}


\def\fonction#1|#2|#3|#4|{\lower 6pt\hbox{
$\matrix{#1 \hfill &\longrightarrow \hfill & #2\hfill\cr
\hfill #3  &\longmapsto \hfill &#4 \hfill\cr}$}}
\newcount\thesection
\newcount\thesubsection
\newcount\thesubsubsection
\newcount\theparagraf
\newcount\thetheo
\newcount\theequ
\global\thesection=0
\global\thesubsection=0
\global\thesubsubsection=0
\global\theparagraf=0
\global\thetheo=0
\global\theequ=0

\font\sevensl= cmsl7
\font\fivesl= cmsl5

\font\tentite = cmbx10 at 16pt
\font\seventite = cmbx7 at 11pt
\font\fivetite = cmbx5 at 8pt
\newfam\titefam
\textfont\titefam = \tentite
\scriptfont\titefam = \seventite
\scriptscriptfont\titefam = \fivetite

\font\sectfont = cmbx10 at 14pt
\font\sectscript = cmbx7 at 10pt
\font\sectsscript = cmbx5 at 7pt
\newfam\sectionfam 
\textfont\sectionfam = \sectfont 
\scriptfont\sectionfam = \sectscript
\scriptscriptfont\sectionfam = \sectsscript
\def\sectionfont{\fam\sectionfam\sectfont}

\font\subsectfont =  cmbx10 at 12pt
\font\subsectscript = cmbx7 at 8pt
\font\subsectsscript = cmbx5 at 6pt
\newfam\subsectionfam
\textfont\subsectionfam = \subsectfont
\scriptfont\subsectionfam = \subsectscript
\scriptscriptfont\subsectionfam = \subsectsscript
\def\subsectionfont{\fam\subsectionfam\subsectfont}

\font\mss=cmss12 scaled \magstep1

\font\hautdepage = cmss8
\font\sf =cmss10
\font\ding = pzdr at 10pt

\font\tenmsb=msbm10
\font\sevenmsb=msbm7
\font\fivemsb=msbm5
\newfam\msbfam
\textfont\msbfam=\tenmsb
\scriptfont\msbfam=\sevenmsb
\scriptscriptfont\msbfam=\fivemsb
\def\Bbb#1{{\fam\msbfam\relax#1}}

\font\teneufrak=eufm10
\font\seveneufrak=eufm7
\font\fiveeufrak=eufm5
\newfam\eufrak
\textfont\eufrak=\teneufrak
\scriptfont\eufrak=\seveneufrak
\scriptscriptfont\eufrak=\fiveeufrak
\def\mathfrak#1{{\fam\eufrak\relax#1}}


\def\begincentered{\par\begingroup
\def \par{\hss\egroup\line\bgroup\hss}\obeylines
\line\bgroup\hss}
\def\endcentered{\hss\egroup\endgroup}

\hsize=12.5cm
\vsize=20cm
\parindent=1cm
\baselineskip=13pt
\hoffset=-0.1cm
\voffset=0.5cm

\long\def\partie#1{\begingroup\sectionfont
        \par\penalty -500
        \vskip 3.25ex plus 1ex minus .2ex
        \skip\afterskip=1.5ex plus .2ex
        \baselineskip=17pt        
        \par
        \def \par{\hss\egroup\line\bgroup\hss}\obeylines
        \line\bgroup\hss
\global\advance\thesection by 1 
\xdef \lastref{\number\thesection}
\global\thesubsection=0 \global\theparagraf=0  \number\thesection\quad#1\hss\egroup\endgroup\par}

\long\def\subpartie#1{\begingroup\subsectionfont%
                          \par\penalty -200
                          \vskip 3.25ex plus 1ex minus .2ex
                          \skip\afterskip=1.5ex plus .2ex
                          \baselineskip=15pt
        \par
        \def \par{\hss\egroup\line\bgroup\hss}\obeylines
        \line\bgroup\hss
\global\advance\thesubsection by 1
\xdef \lastref{\number\thesection.\number\thesubsection}
\global\thesubsubsection=0 \global\theparagraf=0
\number\thesection.\number\thesubsection\quad #1\hss\egroup\endgroup \par}

\def\\{\par}

\long\def\subsubpartie#1{\bgroup\bf
                                \par\penalty -100
                                \vskip 3.25ex plus 1ex minus .2ex
                                \skip\afterskip=1.5ex plus .2ex
                                
\begincentered
\global\advance\thesubsubsection by 1
\xdef \lastref{\number\thesection.\number\thesubsection.\number\thesubsubsection}
\global\theparagraf=0
\number\thesection.\number\thesubsection.%
\number\thesubsubsection\quad #1\endcentered\egroup\par}
\def\paragrafsubsub{%
\global\advance\theparagraf by 1
\xdef \lastref{\number\thesection.\number\thesubsection.\number\thesubsubsection.%
\romannumeral\theparagraf}
{\bf \number\thesection.\number\thesubsection.\number\thesubsubsection.\romannumeral\theparagraf}%
\kern0.2em --- }

\def\paragrafsub{%
\global\advance\theparagraf by 1
\xdef \lastref{\number\thesection.\number\thesubsection.\number\theparagraf}
{\bf \number\thesection.\number\thesubsection.\number\theparagraf}%
\kern0.2em --- }

\def\paragraf{%
\global\advance\theparagraf by 1 
\xdef \lastref{\number\thesection.\number\theparagraf}
{\bf \number\thesection.\number\theparagraf}%
\kern0.2em --- }

\def\paragraphe{%
\par \indent
\ifcase\thesubsection %
  \paragraf
\else
\ifcase\thesubsubsection\paragrafsub %
 \else\paragrafsubsub\fi
\fi}

\def\numeq{\global\advance\theequ by 1%
        \xdef \lastref{(\number\theequ)}%
        \eqno{(\number\theequ)}}
\def\nume{\global\advance\theequ by 1 
        \xdef \lastref {(\number\theequ)}%
        (\number\theequ)}
\def\theo{\global\advance\thetheo by 1%
\xdef\lastref{\number\thetheo}%
\number\thetheo}

\newwrite\fileref
\newread\instream
\newif \ifrefmodif 
\def \defineref#1#2{{\def\next{#1}%
        \expandafter\xdef
           \csname ref = \meaning\next\endcsname{#2}%
        }}

\def \initfileref{
         \openin\instream=\jobname.ref 
         \ifeof\instream \message{ Ahem } 
                \message{**********************************************************}
                \message{******* Le fichier \jobname.ref n'existait pas ***********}
                \message{********** Il faudra absolument recompiler ***************}
                \message{**********************************************************}
                \message{ }
         \else \closein\instream \input \jobname.ref \refmodiffalse 
         \fi 
        \ifarxiv\else\immediate\openout \fileref=\jobname.ref\fi
        \global\let\initfileref=\relax} 

\def \label#1{\ifarxiv\initfileref\else{%
        \toks0={#1}\wlog{REF \the\toks0= \lastref}
        \initfileref
        \immediate\write\fileref{\noexpand\defineref%
                {\the\toks0}{\lastref}}%
        \def\next{#1}%
        \expandafter\ifx 
          \csname ref = \meaning\next\endcsname\lastref
        \else \global\refmodiftrue  \message{ }\message{ Attention }
        \message{reference {\the\toks0 = \lastref }  modifiee ou redefinie} \message{ }\fi 
        \defineref{#1}{\lastref}%
        \ifrefvis\ifhmode\raise 6pt \vbox{\hbox to 0pt{\hss\fivebf [\the\toks0]\hss}}
        \else\llap{\fivebf [\the\toks0]}\fi\fi}\fi}

\def \ref#1{{\initfileref \def\next{#1}%
        \csname ref = \meaning\next\endcsname}}

\outer\def \bye{
        \closeout \fileref \vfill \supereject
        \ifrefmodif \erreurmodif \fi
        \end}

\def \erreurmodif{
        \message{ }
        \message{**********************************************************}
        \message{BEWARE, some references have been modified.}%
        \message{Please compile the tex file again to get references right,}
        \message{if this message appear again, then a reference must have}
        \message{been define at least two times.}        
        \message{**********************************************************}
        \message{ }}


\def \ignorepar{\afterassignment\ignoreparaux \let\next=}
\def \ignoreparaux{\ifx\next\par \let\next\ignorepar \fi \next}

\def\bibliography{\bgroup\sectionfont
        \par\penalty -500
        \vskip 3.25ex plus 1ex minus .2ex
        \skip\afterskip=1.5ex plus .2ex
\begincentered References \endcentered\egroup\par
\tabskip5pt minus 1pt%
\noindent\halign to \hsize \bgroup##\hfil&##\hfil\cr}
\def\bibliend{\egroup}

\def\bibitem[#1]#2{\xdef\lastref{[#1]}\label{#2}\hskip-8pt[#1]&%
\vtop\bgroup\hsize10.6cm\noindent\ignorespaces}
\def\finitem{\medskip \egroup \cr}

\def\newblock{\hskip .11em plus .33em minus .07em}
\def\cite#1{\ref{#1}}

\def\em{\sl}

\def\proof{\noindent{\bf Proof. }}
\def\qed{\ \-\hfill {\ding \char 111} \par}

\def\C{{\Bbb C}}
\def\T{{\Bbb T}}

\def\R{{\Bbb R}}
\def\Z{{\Bbb Z}}
\def\N{{\Bbb N}}
\def\frac#1#2{{#1\over#2}}
\def\tend{\mathop{\longrightarrow}}

\def \scal#1{\prodscal #1>} 
\def \prodscal#1,#2>{\langle #1,#2 \rangle}

\def\n#1{||#1||}
\def\bign#1{\bigl\|#1 \bigr\|}

\def\Biggv#1{\Biggl|#1 \Biggr|}

%
%
